\documentclass[]{interact}
\usepackage{latexsym}
\usepackage{anyfontsize} %the T&F interact.cls file has some non-integer font sizes (like 5.5) so this fixes that warning

\usepackage{xcolor}
\usepackage{array}
\usepackage{graphicx}
\usepackage{tabularx}% http://ctan.org/pkg/tabularx
\usepackage{booktabs}% http://ctan.org/pkg/booktabs
\usepackage{verbatim}
\usepackage{comment}
\usepackage{amssymb}
\usepackage{mathrsfs} % some font support
\usepackage[export]{adjustbox}
\usepackage{hyperref}
\usepackage{makecell}

\usepackage{comment}
 
\theoremstyle{plain}
\numberwithin{equation}{section} %change this to make globally numbered environments
\newtheorem{thm}{Theorem}[section]
\newtheorem{theorem}[thm]{Theorem}
\newtheorem{lemma}[thm]{Lemma}

\newtheorem{definition}[thm]{Definition}

\newtheorem{corollary}[thm]{Corollary}
\newtheorem{algorithm}[thm]{Algorithm}
\numberwithin{equation}{section}

\DeclareMathOperator{\tree}{Tree}
\DeclareMathOperator{\shape}{Shape}
\newcommand{\ttt}[1]{%
\ensuremath{\mathfrak{#1}%
}}
 
\DeclareMathOperator{\tcount}{F}
\DeclareMathOperator{\generatingfunction}{\mathcal{F}}
\DeclareMathOperator{\rodcount}{C}
\DeclareMathOperator{\rodgeneratingfunction}{\mathcal{C}}
\DeclareMathOperator{\trains}{Trains}

\DeclareMathOperator{\sign}{sign}
\DeclareMathOperator{\length}{length}
\DeclareMathOperator{\rodcolor}{color}

\DeclareMathOperator{\discrepancySeq}{D}

\newcommand{\RQS}{\ensuremath{R \trailright{Q} S}}

\newcommand{\anti}[1]{%
\ensuremath{\overline{#1}}%
}
\newcommand{\rod}[1]{%
\ensuremath{|#1|}%
}
\newcommand{\oddss}[1]{%
\ensuremath{#1^o}%
}

\newcommand{\seq}[1]{%
\ensuremath{\mathscr{#1}%
}}
\newcommand{\aaa}{%
\ensuremath{\alpha}%
}

\newcommand{\emptyrods}{[\phantom{a}]}
\newcommand{\emptytrain}{\Lambda}
\newcommand{\trailright}[1]{\xrightarrow{\hphantom{a}#1\hphantom{a}}}
\newcommand{\inner}{\text{Inner}}
\newcommand{\leaves}{\text{Leaves}}
\newcommand{\size}[1]{\ensuremath{\Vert #1 \Vert}}

\begin{document}

\title{Recursions, Trains, Trees, and Combinatorial Rod Set Algebra}

\author{
\name{Ethan D. Bolker\textsuperscript{a},
Debra K. Borkovitz\textsuperscript{b}
\thanks{CONTACT  Debra K. Borkovitz Email: dbork@bu.edu} 
and
  Katelyn Lee
}
\affil{
\textsuperscript{a}Department of Mathematics, University of Massachusetts, Boston,
\textsuperscript{b}Department of Mathematics and Statistics, Boston University
}
}
\date{\today}

\maketitle

{\bf Article type}: mathematical outreach 

\bigskip

\begin{abstract}
  We explore a physical model of ordered sums of integers as  trains  of rods. The trains for a fixed, possibly infinite, set of rods naturally correspond to nodes in a tree;  relations among finite linear recursions encoded in the subtrees define  algebraic operations on sets of rods. We use this algebra to prove classic identities for recursively defined sequences, to show that some Lucas sequences are divisibility sequences, to characterize two-term linear Fibonacci identities, and to find the cyclotomic polynomial factors of Borwein trinomials. We complement abstractions with lots of examples.
\end{abstract}


\section*{Introduction}

\begin{figure}[h]
  \centering
\includegraphics[width=0.4\textwidth]{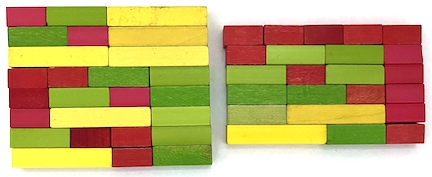}
\caption{The $14$  trains of length $10$ made from  Cuisenaire\texttrademark{} rods of lengths
  $2$ (red), $3$ (green), and $5$ (yellow). Each row is a separate train. }  
  \label{fig:Cuisenaire}
 \end{figure}

Children use
Cuisenaire\texttrademark{} rods~\cite{cuisenaire} to learn how numbers
fit together. A common activity is to make trains consisting of a sequence of rods placed end-to-end, and to count how many trains are possible under various restrictions. 
Figure~\ref{fig:Cuisenaire} shows all possible trains of length $10$ made from 
rods of lengths 
     $2$,  $3$, and $5$. The trains  model \emph{compositions}: sequences of positive integers with a given sum.

The Fibonacci numbers count the ways to build a train of length $n$
from rods of lengths $1$ and $2$; the recursion
\begin{equation*}
\tcount(n) = \tcount(n-1) + \tcount(n-2)
\end{equation*}
sorts the trains of length $n$ depending on the length of the first (or last)
rod. In this model for the Fibonacci sequence, indexing starts at $1$ with  $\tcount(1)=1$ and $\tcount(2) = 2$.

Playing with the recursion, we see that
\begin{align*}\tag{1} %\label{eq:fib134}
  \tcount(n) &=\tcount(n-1)+\tcount(n-2) \notag \\
  &= \tcount(n-1)+(\tcount(n-2-1)+\tcount(n-2-2)) \notag \\
  &= \tcount(n-1)+\tcount(n-3)+\tcount(n-4) .
\end{align*}

It follows that for $n>2$, the sequence $\tcount(n)$
also satisfies the
recursion determined by rods of length
 $1$, $3$, and $4$.

\begin{figure}[h!]
    \centering
    \includegraphics[width=0.65\linewidth]{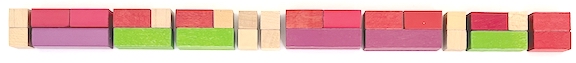}
        \caption{Converting the top train of $1$'s and $2$'s into the bottom train of $1$'s, $3$'s, and $4$'s, with a $2$ at the end. }
    \label{fig:longrodphoto}
\end{figure}

Soon we will use our rod model to study this observation about the Fibonacci sequence by converting
trains composed of rods of lengths $1$ and $2$ into trains composed of rods of lengths $1$, $3$, and $4$. When a $2$ is followed by a $1$ we can glue them together to make a $3$, and when a $2$  is followed by a $2$ we glue to make a $4$. This construction may end with an extra rod of length $2$ at the end of the train.   Figure~\ref{fig:longrodphoto} shows an example.

 Kids can't play with infinitely many rods, but we can.
 If we repeatedly replace the rods of even length $2k$ by rods of
lengths $2k+1$ and $2k+2$ we are led to the recursion defined by rods of lengths
$[1,3,5,7, \ldots]$. We will see how to define the Fibonaccis so they satisfy that recursion when $n$ is odd.

Kids can identify Cuisenaire\texttrademark{} rods by  color or by length. But
writing
\begin{align*}\tag{2}%\label{eq:fib223}
  \tcount(n) &=\tcount(n-1)+\tcount(n-2) \notag \\
  &= (\tcount(n-1-1)+\tcount(n-1-2))+\tcount(n-2) \\ \notag 
  &= 2\tcount(n-2) + \tcount(n-3) \notag,
\end{align*}
  shows that when $n > 3$, the Fibonaccis also satisfy the
  recursion for 
   trains built from two different kinds of rods of length $2$
  along with 
  rods of length $3$. To preserve our physical model we will
  allow rods with the same length
  but different colors.   

  Since
  \begin{equation*}
    \tcount(n-1) = \tcount(n-2) + \tcount(n-3),
  \end{equation*}
  we know that
  \begin{equation*}
    \tcount(n-2) = \tcount(n-1) - \tcount(n-3),
  \end{equation*}
  so
\begin{align*}\tag{3}%\label{eq:fib11anti3}
  \tcount(n) &=\tcount(n-1)+\tcount(n-2) \notag \\
  &= \tcount(n-1) + (\tcount(n-1) -\tcount(n-3)) \notag \\
  &= 2\tcount(n-1)-\tcount(n-3),
\end{align*}
when $n > 3$,
which suggests that we will also need a way to define ``negative" rods (which will not be rods of negative length).

In the first five sections we'll formalize this rod train model for relationships between recursions. Section 1 has careful definitions for rods and rod trains. Sections 2 and 3 introduce a tree structure for rod trains that leads to a natural way to understand the recursion manipulations in this Introduction.
In Sections 4 and 5 we introduce algebraic operations on rod sets. 

In Section 6 we pause to connect our work to generating functions; we see there that rod set algebra 
provides a combinatorial interpretation for polynomial multiplication. In Section 7 we define duality, the last operation on rod sets, and use it to prove many identities.  The next six sections cover applications, 
often involving  small finite rod sets. We study periodic sequences and connect patterns of train counts to rod set operations, which leads to more Fibonacci identities, a new proof that some Lucas sequences are divisibility sequences, and a new way to understand the cyclotomic factors of Borwein trinomials. We end with further avenues to explore. 

To the best of our knowledge our techniques are new.

\section{Setting the Stage}\label{sec:setting}

We start with careful definitions and some elementary observations about rods, trains and sets of rods and  trains.

 \begin{definition}\label{def:rod}
  A \emph{rod} is a triple 
  \begin{equation*}
      r = (\length(r), \rodcolor(r), \sign(r)),
  \end{equation*}
where $\length(r)$ is a
  positive integer,  $\rodcolor(r)$ is an arbitrary tag we think of as a
  \emph{color}, and 
 $\sign(r) = \pm 1$.
 
 For rod $r$, we write $\anti{r}$ for the rod with the same length and color and opposite sign.
 A rod with negative sign is an \emph{antirod}.
 Sometimes ``rod'' will mean a rod of either
sign and sometimes just a rod with positive sign. The context should make clear when the sign is important.

\end{definition}

In  explicit examples we abbreviate rods by writing just their lengths, so $2$ represents the rod $(2, c, 1)$ and  $\anti{2}$ represents the rod $(2,c,-1)$. When we need to distinguish rods by color we write $2_c = (2, c, 1)$ and $\anti{2}_c = (2, c, -1)$. Often we  use the length of a rod $r$ in an arithmetic expression such as $n-\length(r)$. When no confusion can result we will omit the $\length$ function and write simply $n-r$. Formally,  if $n$ is an integer and rod $r = (m,c, \pm 1)$ then $n-r = n-m$.

\begin{definition}\label{def:rodtrain}
A \emph{rod train}, or  simply a \emph{train},  is a finite sequence of rods.
The \emph{empty train} is the empty sequence $\emptytrain$. We write trains as concatenations. 

For example,
    $\anti{1}\, \anti{2}\, \anti{1}$, $1\anti{2}1$ and  $1_c1_d\anti{2}$ are trains. In the first two the implicit assumption is that the two rods of length 1 have the same color.

For the nonempty train $\tau$ we write $\rod{\tau}$ for
the rod with length the sum of the lengths of the rods in $\tau$, sign the product of the signs, and color that distinguishes it from other rods of the same length and sign in its context. The colors of the rods in $\tau$ have no effect on $\rod{\tau}$.
Think of this operation as
gluing together the rods in $\tau$. 
Since by definition rods have positive length, there is no rod
$\rod{\emptytrain}$.
For example, 
$\rod{\anti{1}\,\anti{2}\,\anti{1}}=\rod{1\anti{2}1} = \rod{1_c1_d\anti{2}}  = \anti{4}$.

We define the \emph{length} and \emph{sign} of train $\tau$ as the length and sign of rod $\rod{\tau}$. Thus a train is \emph{positive}
  when it contains an even number of antirods and \emph{negative}
  otherwise. 

We concatenate trains $\sigma$ and $\tau$ as either $\sigma\tau$ or $\sigma.\tau$; we use the dot when we want it for emphasis. 
In general $\sigma\tau \ne \tau\sigma$, but $\rod{\sigma\tau} = \rod{\tau\sigma}$ always.
\end{definition}

\begin{definition}
A \emph{rod set} $R$ is a set of rods with 
finitely many rods of
each length. We write a rod set as a square
bracketed list of rods, usually in non-decreasing length
order. Duplicate  entries implicitly specify rods of
different colors.
When there are multiple rods of the same length we may use exponents to count them, so
\begin{equation*}
    [1, \anti{1},\anti{1},2,2,2] = [1, \anti{1}^2,2^3].
\end{equation*}
We write $\anti{R}$ for the rod set whose rods reverse
the signs of the rods in $R$.

The \emph{concatenation} of  rod sets $R$ and $S$ is the rod set
\begin{equation*}
    R.S = RS = SR = \{   \rod{jk} \mid j \in R, k \in S \}, 
\end{equation*}
so for example
\begin{equation*}
    [1, \anti{1},\anti{1},2,2,2].[1,5] =  
    [2, \anti{2}^2,3^3,6,\anti{6}^2,7^3].
\end{equation*}

 Let $R^+$ and  $R^-$ be the sets of rods and antirods in $R$ and 
 $R^+(n)$ and $R^-(n)$ the finite sets of positive and negative rods of length $n$. 
  Writing $\#S$ for the cardinality of $S$, define
  the \emph{net rod count} $\rodcount(n,R)$ for $n > 0$ as:
  \begin{equation*}
      \rodcount(n,R) =  \#{R^+(n)}- \#R^-(n).
  \end{equation*}
 \end{definition}
 The net rod count $C(n,R)$ can be positive, negative, or zero.
 For example, for the rod set $R = [1,1, \anti{1},\anti{2},3,\anti{3}]$, we have
 $\rodcount(1,R) = 1$, $\rodcount(2,R) = -1$, and all the other net rod counts are $0$.

 Since there are no rods with length $0$ we could consistently define
 $\rodcount(0,R) = 0$. For reasons that will become clear later on, we choose instead to leave it undefined.

In    rod set unions we assume that colors are all distinct, so
\begin{equation*}
    [2, \anti{2},3] \cup [2,2,4] = [2^3,\anti{2},3,4].
\end{equation*}
In intersections and set differences we assume colors match as much as possible, so
\begin{equation*}
    [2, 2,\anti{2},3] \cap [2,2,4] = [2,2]  
\end{equation*}
and
\begin{equation*}
    [2,2,\anti{2},3] \setminus [2] =   [2,\anti{2},3] .
\end{equation*}

\begin{definition}
    A \emph{rod train set}, or simply a \emph{train set}, is a set $T$ of rod trains with only finitely many trains of each length $n$. The \emph{net train count} $\rodcount(n,T)$ is the difference between the number of positive and negative trains of length $n$. 

    For train set $T$ we write
    \begin{equation*}
        \rod{T} = \{ \rod{\tau} \mid \tau \in T\}
    \end{equation*}
for the set of rods built from the nonempty trains in $T$.
Thus $\rodcount(n,T) = \rodcount(n, \rod{T})$ when $n > 0$.

     For example, $  
   T = \{
    11\anti{2}1\anti{2}1   ,
     11\anti{2}1\anti{2}11  ,
     1\anti{2}1  \}
$
is a train set built from the rods in $[1,\anti{2}]$. Here
$\rod{T} = [8,9,\anti{4}]$ and 
$\rodcount(4,T) = -1$.

  The  infinite set of all trains built from the rods in rod set $R$ is   $\trains(R)$: all the finite sequences of rods in $R$. We refer to these as $R$-trains.
\end{definition}

Now we are ready to define the integer sequences whose recursive properties we will explore.

\begin{definition}\label{def:traincounts}
  For rod set $R$, the \emph{net train count sequence} is
  \begin{equation*}
       \tcount(n,R) =  \rodcount(n, \trains(R)) =
       \begin{cases}
               \rodcount(n, \rod{\trains(R)})  & n > 0 \\
               1 & n = 0 \\
               0 & n < 0.
       \end{cases}
\end{equation*}
We usually think of this sequence as starting at $n=0$. The $0$ values for $n < 0$ are there to make some algebra cleaner.

When the rod set $R$ is clear from context we may shorten expressions by writing $\tcount(n)$ instead of 
$\tcount(n,R)$.
 \end{definition}

 The traditional function identifier, ``$\tcount$", for the function that counts trains is a nod to the start of our study. It's short for \emph{Fibonacci}: the sequence of train counts for rods of lengths $1$ and $2$. 

Table~\ref{table:1anti2} displays the positive and negative trains of length at
most $5$ for the rod set $[1,\anti{2}]$. Ignoring signs, the total number of trains
of each length is, of 
course, the corresponding Fibonacci number.  The last column is the net train count.

 \begin{table}[h!]
     \centering
     \caption{Rod trains for $[1, \anti{2}]$}.
     
     \begin{tabular}{|r|r|r|r|r|}
\hline 
     \thead{Length} & \thead{Positive} & \thead{Negative} & \thead{Train \\ Count} & \thead{$\tcount(n,[1,\anti{2}])$} \\
     \hline
     0 &  $\emptytrain$  &  & 1 & 1\\
     1 & 1 &  & 1 & 1\\
     2 & 11                    & \anti{2}                         & 2 &  0 \\
     3 & 111                   &1\anti{2},\anti{2}1               &3  & -1 \\
     4 & 1111,\anti{2}\anti{2} & 11\anti{2},1\anti{2}1,\anti{2}11 & 5
     & -1 \\
     5 & 11111,1\anti{2}\anti{2},\anti{2}1\anti{2},\anti{2}\anti{2}1 & 111\anti{2},11\anti{2}1,1\anti{2}11,\anti{2}111 & 8
     & 0 \\     
     \hline
     \end{tabular}
     \label{table:1anti2}
  \end{table}

  \begin{theorem}\label{thm:traincountsarerecursions}
  For rod set $R$, the net train count sequence $\tcount(n,R)$ satisfies the recursion
\begin{equation}\label{eq:recursion}
  \tcount(n,R) =     \sum_{ k \in R} (\sign k)\tcount(n-k, R),
\end{equation}
with initial conditions
\begin{equation}\label{eq:initialconditions}
\tcount(n,R) =
\begin{cases}
  1 \text{ if } n = 0 \\
  0 \text{ if } n < 0 \ .
\end{cases} 
\end{equation}
\end{theorem}

\begin{proof}
Let $ \tcount^+(n,R)$ be the number of positive trains of length $n$ in  $\trains(R) $  and
$ \tcount^-(n,R)$ the number of negative trains.
Then
     \begin{equation*}
     \tcount(n,R) = \tcount^+(n,R) - \tcount^-(n,R).      
     \end{equation*}

 Let $k$ be the last rod in a train $\tau$ of length $n>0$. Removing that rod from $\tau$ leaves a train $\tau'$ of length $n - \length(k)$, which we write as $n-k$. If $k$ is an antirod the signs of $\tau$ and $\tau'$ are opposite, otherwise they are the same.
 Then

\begin{equation*}
    F^+(n,R)=\sum_{k \in R^+}F^+(n-k,R) + \sum_{k \in R^-}F^-(n-k,R), 
\end{equation*}
and

\begin{equation*}
    F^-(n,R)=\sum_{k \in R^-}F^+(n-k,R) + \sum_{k \in R^+}F^-(n-k,R). 
\end{equation*}
Note that if $k>n$ the initial conditions for negative $n$ ensure that no extra terms are added. 

Thus
\begin{align}
    \tcount (n,R) &=\tcount^+(n,R) -  \tcount^-(n,R)  \notag \\
       &= \sum_{k \in R^+} \Bigl( \tcount^+(n-k, R) - \tcount^-(n-k, R)\Bigr) \notag \\
        & \ \ \ \ \ 
       - \sum_{k \in R^-} \Bigl(\tcount^+(n-k, R) - \tcount^-(n-k, R)\Bigr)  \label{eq:signchange} \\
   &=  \sum_{k \in R} (\sign k)\Bigl (( \tcount^+(n-k, R) - \tcount^-(n-k, R)\Bigr) .\notag
\end{align}
\end{proof}

The initial condition $\tcount(0,R) = 1$ 
guarantees that when $k \in R$ the one-rod train $k$
contributes $\pm 1$ to the net train count.  
It is consistent with the fact that the
empty train contains an even number of antirods. At the end of Section~\ref{sec: algebra} we show how to use rod trains to model recursive sequences with different initial conditions.

In Benjamin and Quinn's \emph{Proofs that Really Count} \cite{benjamin},
rod sets with these default initial conditions correspond to
\emph{weighted unphased tilings} of a board. 
Because we allow infinite rod sets and are particularly interested in the relations among rod sets, we  focus on the rods themselves rather than introducing a board for them to tile. Their system of weights accommodates any complex number, whereas we only need weights of $\pm 1$. Hopkins \cite{hopkins} also uses Cuisenaire rods in his many engaging classroom activities, some of which overlap with our work. We branch off when we include antirods and study the connections among rod sets.

Here are some familiar sequences described as net train counts, starting at $n=0$:

 \begin{itemize}
     \item $\tcount(n, [1,2])= 1,1,2,3,5,8,\ldots$,  the 
       Fibonacci numbers.      
\item $\tcount(n, [\anti{1},2])= 1,-1,2,-3,5,-8,\ldots$,  the 
      alternating Fibonacci numbers.   
     \item $\tcount(n, [1,3,5,7,\ldots])= 1,1,1,2,3,5,8,\ldots$, the
       Fibonacci sequence with an extra initial $1$.
     \item $\tcount(n,[1]) = 1,1,1,1,\ldots$, the constant sequence of all $1$'s.
     \item $\tcount(n,[1,1]) = 1,  2, 4, 8, 16, \ldots$, the powers of $2$.
     \item $\tcount(n,[2,3]) =  1, 0, 1, 1, 1, 2, 2, 3, 4, 5, 7, 9, 12, \ldots $, a \emph{Padovan
   sequence} \cite[p 89]{hopkins} \cite{padovan}.
     \item $\tcount(n,[1,3]) = 1, 1, 1, 2, 3, 4, 6, 9, 13, 19, \ldots$, a \emph{Narayana cow number
     sequence} \cite{narayana-oeis}.
          \item $\tcount(n,[1,\anti{2}]) = 1, 1,0,-1,-1,0, 1, \ldots$, repeating
       with period  $6$.
       \item $\tcount(n,[\anti{1},\anti{2}]) = 1, -1,0,1,-1,0, 1, \ldots$, repeating
       with period  $3$.
       \item 
       $\tcount(n,\emptyrods)= \tcount(n,[1,\anti{1}]) = 1, 0,  0, 0, \ldots$ .
 \end{itemize}

The next theorem shows that every integer sequence is the net train count sequence for some rod set.

\begin{theorem}\label{thm:recursionsarerodsets}
    For every infinite sequence of integers
  $ a_1, a_2, \ldots $
    there is a rod set $R$ such that  $\tcount(n,R)= a_n$   for all $n \ge 1$.
\end{theorem}
\begin{proof}
We construct $R$ inductively. By definition, $\tcount(0,R) = 1$.
To make $\tcount(1,R) = a_1$, the set $R$ clearly needs $a_1$ differently colored rods of length
$1$ (these will be antirods if $a_1 < 0$). Once the number of rods of length $k$ is known for all $k < n$, we can use the recursion in Equation~\ref{eq:initialconditions} to compute the number of trains of length $n$ composed of more than one rod. Then we add the correct number of differently colored
rods or antirods of length $n$ so that these new one-rod trains make
  $\tcount(n,R) = a_n$.
\end{proof}

For the train counts corresponding to the Catalan numbers indexed from $0$:
\begin{equation*}
\tcount(n,R) = 
1, 1, 2, 5, 14, 42, 132, 429, \ldots,
\end{equation*}
 the corresponding rod set is
\begin{equation*}
   R =   [1, 2, 3^2, 4^5, 5^{14}, 6^{42}, 7^{132}, 8^{429}, \ldots ].
\end{equation*}
The number of rods of each length is the Catalan number  with index shifted by $1$ from the train counts. The recursion corresponds to the well-known convolution that the Catalans satisfy.

\section{Trees}\label{sec:trees}

Here we explore a natural way to organize
$\trains(R)$ 
in a tree\footnote{
Sorry for the awkward   metaphor mixing trees and trains.} 
that lets us systematically study the way we combined rods in Figure~\ref{fig:longrodphoto}.

\begin{definition}\label{def:Rtree}
 Let $R$ be a rod set.   $\tree(R)$ is the rooted tree in which each node has one child for each rod in $R$. We label the edge from parent to child by that rod. Then we label each node by the rod train that is the path from the root to that node. The root's label is the empty train $\emptytrain$. 
 The children of a node labeled by the train $\tau$ are labeled by the  trains $\tau.R$ constructed by appending a rod $k \in R$ to $\tau$.
The nodes at the $m$th level correspond to trains containing $m$ rods.

When $\tau \in \trains(R)$ and no confusion will result, we may write ``node $\tau$" for the corresponding node in $\tree(R)$ instead of the more formal ``node with label $\tau$." 

We decorate node $\tau$ by adding $\rod{\tau}$ as a second label. 
In practice, these second labels do not identify nodes, because it is cumbersome and usually unnecessary to assign explicit different colors to rods of the same length in the rod set
     \begin{equation}\label{eq:geometrictrains}
        \rod{\trains(R)} =   R \cup RR \cup RRR \cup \cdots .
    \end{equation}
\end{definition}

In the language of theoretical computer science, $R$ is an \emph{alphabet} and $\trains(R)$ are the \emph{words} or \emph{strings} in the \emph{language} for that alphabet. 
The tree structure depends only on the number of letters in the alphabet (which may be infinite). When $R$ contains just two rods its tree is the complete infinite binary tree.

\begin{figure}
    \centering
    \includegraphics[width=\textwidth]{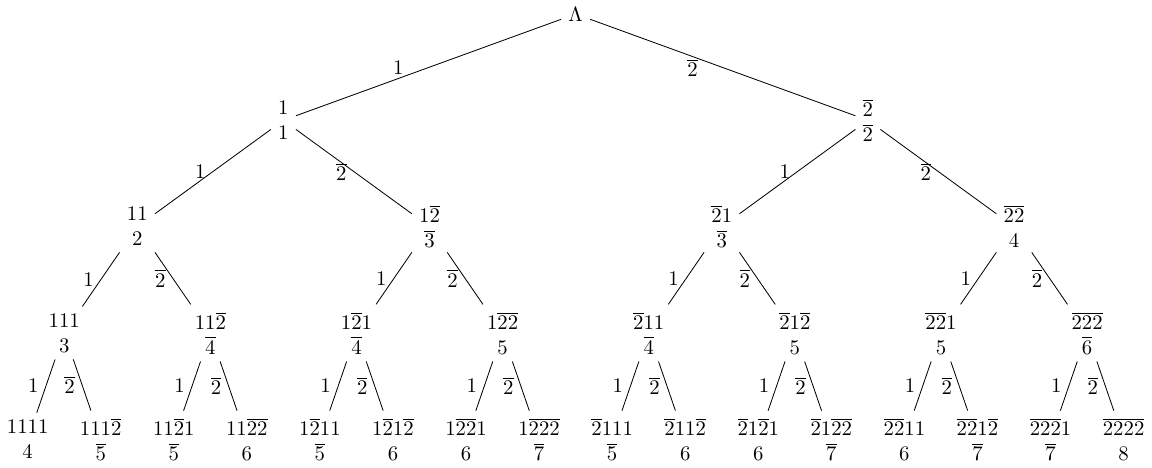}
    \caption{The first four levels of $\tree([1,\anti{2}])$}
    \label{fig:1anti2trains}
\end{figure}

Figure~\ref{fig:1anti2trains} shows the first four levels of 
$\tree([1,\anti{2}])$.  The Fibonacci numbers count the number of node labels of each length. We can use this figure to confirm the entries in Table~\ref{table:1anti2}.

 Figure~\ref{fig:longrodphoto} 
in the Introduction corresponds to 
replacing the node   $\tau  = 2$   at the first level of $\tree([1,2])$ by its children $21$ and $22$. We write this operation on nodes as
\begin{equation}\label{eq:12expand2}
\{1,2\}\trailright{2}\{1,21,22\} .
\end{equation}

We use Equation~\ref{eq:12expand2} to replace some pairs of rods in an element of $\trains([1,2])$ with an element of $\trains([1,3,4])$, possibly with an extra rod of length $2$ at the end. For the example in Figure~\ref{fig:longrodphoto} the parsing is
\begin{equation}\label{eq:12expand2trains}
    \tau = 12221211122221212=1.22.21.21.1.1.22.22.1.21.2,
\end{equation}
which maps $\tau$
to the element $1433114413$ of $\trains([1,3,4])$ with an appended rod of length $2$.

We can iterate this kind of tree node replacement. For example, replacing node $1$ next and then node $21$ leads to the expansion chain we write as
\begin{equation} \label{eq:trailwithtrains}
\{1,2\}\trailright{2}\{1,21,22\} \trailright{1} \{11,12,21,22\}
\trailright{21} \{11,12,211,212,22\} .
\end{equation}

 The next definition formalizes
this kind of expansion.

\begin{definition}\label{def:treeexpansion}
A rooted subtree $\ttt{T}$   of $\tree(R)$ is \emph{an expansion subtree} when
\begin{itemize}
    \item It contains all the children of the root.
    \item If it contains a child of some node, then it contains all the children of that node.
\end{itemize}
For an expansion subtree $\ttt{T}$,  $\inner(\ttt{T})$ is the set of nodes (other than the root) of $\ttt{T}$ that have at least one child (hence all children) in $\ttt{T}$; we call these the \emph{internal nodes}. $\leaves(\ttt{T})$ is the set of nodes of $\ttt{T}$ that have no children in $\ttt{T}$.

We write the expansion corresponding to expansion subtree $\ttt{T}$ as
\begin{equation*}
    R \trailright{\inner(\ttt{T})} \leaves(\ttt{T}) .
\end{equation*}
\end{definition}

The expansion subtree $\ttt{T}$ whose only nodes are the children of the root has $\leaves(\ttt{T}) = R$ and empty $\inner(\ttt{T})$.  $\tree(R)$ itself is an expansion tree with no leaves. All its nodes except the root are internal.

\begin{figure}  
    \centering   
          \includegraphics[ width= 0.26\textwidth,valign=t]{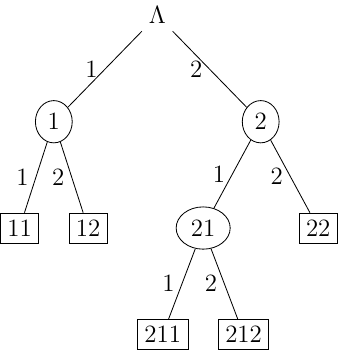}    
     \phantom{aa}
    \includegraphics[height=2in,valign=t]{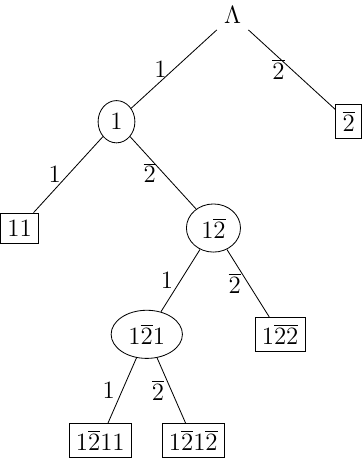}
     \phantom{aa} 
     \includegraphics[height=2in,valign=t]{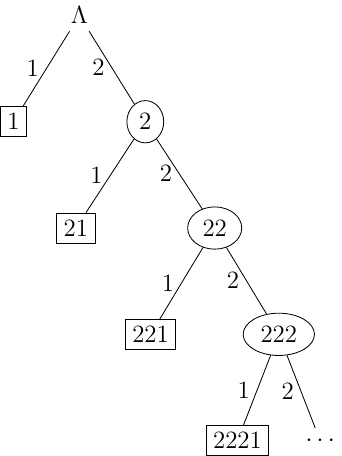}
  \caption{ Expansions of $[1,2]$,    $[1,\anti{2}]$, and $[1,2]$,
  with internal nodes circled, leaves boxed. }
  \label{fig:qs-expansion}    
\end{figure}

Figure~\ref{fig:qs-expansion} shows three  expansions; in each figure the internal nodes are circled and the leaves boxed:
\begin{align*}
  \{1,2\} & \trailright{\{1, 2, 21\} } \{11, 12, 22, 211, 212\},\\
    \{1,\anti{2}\} & \trailright{\{1 ,1\anti{2} ,1\anti{2}1\}} 
     \{\anti{2}, 11,  1 \anti{2} \anti{2},  1 \anti{2}11, 1 \anti{2}1\anti{2} \},\\
    \{1,2\} & \trailright{\{2, 22, 222, \ldots\}}
    \{1, 21, 221, 2221, \ldots \}.
\end{align*}

\section {Discrepancies}

In the Introduction we saw that the Fibonacci train counts  $\tcount(n,[1,2])$ satisfy the recursion  for the rod set $[1,3,4]$  except when $n=2$, which was the single rod we replaced in Expansion~\ref{eq:12expand2}.

In this section we show that the internal nodes in an expansion subtree correspond to the places where the original net train counts fail to satisfy the recursion determined by the expansion.
  
 We begin with a definition that refers only to rod sets and 
 %recursions, 
 their corresponding recurrence relations,
 not to subtrees.

\begin{definition} \label{def:discrepancy}
    For rod sets $R$ and $S$ the \emph{discrepancy} at $n>0$ is the difference
    \begin{equation*}
    \discrepancySeq(n,R,S) = \tcount(n,R)- \sum_{k \in S} (\sign k)\tcount(n-k,R).
\end{equation*}
\end{definition}

Thus $\discrepancySeq(n,R,S)$ measures by how much the net train counts for $R$ fail to satisfy the recursion for $S$ at $n$. We look at some examples to get a feel for discrepancies and to motivate the theorem that will characterize them.

 First, 
 $\discrepancySeq(n,[1,2], [1,3,4])$ is $1$ when $n=2$ and is $0$ everywhere else.

 Clearly $\discrepancySeq(n,R,R)=0$ for all $n>0$.

Earlier we asserted without proof that the net train counts
\begin{equation*}
    \tcount(n,[1, \anti{2}]) = 1,1,0,-1,-1,0, \ldots  
\end{equation*}
repeat with period  $6$. A sequence with period $6$ satisfies the recursion defined by the rod set $[6]$, 
so to prove that assertion we will show that for $n>0$ the discrepancies satisfy
\begin{equation}\label{eq:1anti2-6}
    \discrepancySeq(n,[1, \anti{2}],[6]) = 
    \begin{cases}
        1 & n = 1 \\
        -1 & n = 3,4 \\
        0 & \text{otherwise}.
    \end{cases}
\end{equation}
This will follow from the next theorem, which will also imply the following two 
  Fibonacci identities for the sequence $\tcount(n) = \tcount(n, [1,2]) = 1,1,2,3,5,8,13,\ldots$.
  \footnote{
  The traditional indexing for this sequence starts with $F_0=0$, then
  $F_1 = F_2 = 1$, thus shifted one from ours.
  }

\begin{align}\label{eq:Fibidentity1}
    1 + 2 + 5 + 13 +  34 +  \ \cdots   \   +\tcount(2n)    &= \tcount(2n+1) ,\\
    1 + 3 + 8 + 21 + \cdots + \tcount(2n-1) &= \tcount(2n) -1 . \label{eq:Fibidentity2}
\end{align}

Since the recursion for the rod set $ \text{Odds} = [1,3,5,\ldots]$ is
\begin{equation*}
    \tcount(n, \text{Odds}) =   \tcount(n-1, \text{Odds})  + \tcount(n-3, \text{Odds}) + \cdots ,
\end{equation*}
the above identities are equivalent to
\begin{equation*}
    \discrepancySeq(n, [1,2], \text{Odds})  = 
    \begin{cases}
        0 &\text{if $n$ is odd} \\
      1 &\text{if $n$ is even},
    \end{cases}
\end{equation*}
 which will follow from the next theorem. 

\begin{lemma}\label{lemma:parseRtrains}(Parsing)
    Suppose  $\ttt{T}$ is an expansion subtree of  $\tree(R)$.   Then 
    we can write every nonempty $R$-train $\rho$  uniquely as a concatenation
    \begin{equation*}
        \rho =    \sigma_1.\sigma_2. \cdots .\sigma_m.\tau
    \end{equation*}
    where each $\sigma_i \in \leaves(\ttt{T})$ and   $\tau \in \inner(\ttt{T})$. 
    In this expression $m$ may be $0$ (that is, no $\sigma_i$s) and
    $\tau$ may be $\emptytrain$, but not both. 
    
    Conversely, every such expression is an $R$-train.
\end{lemma}
\begin{proof}
The argument generalizes the example in Equation~\ref{eq:12expand2trains}. Given $\rho \in \trains(R)$, we follow $\rho$ from the root of $\tree(R)$  until we encounter a node $\sigma_1 \in \leaves(\ttt{T})$. We continue inductively with the rest of $\rho$ to find $\sigma_2,  \ldots \sigma_m$. The construction ends when there are no rods left, or too few left to reach a node in $\leaves(\ttt{T})$. Then the 
trailing train ends at a node $\tau \in \inner(\ttt{T})$.

The converse is the simple observation that each $\sigma_i$ is just a sequence of rods from $R$.
\end{proof}

Using the second expansion tree in Figure~\ref{fig:qs-expansion} we can see the parsings
\begin{equation*}
    \begin{array}{llll}
    11\anti{2}1\anti{2}1&=  11.\anti{2}.1\anti{2}1   &\tau = 1\anti{2}1 &m = 2,\\
     11\anti{2}1\anti{2}11&=  11.\anti{2}.1\anti{2}11 &\tau = \emptytrain  &m = 3,\\   
     1\anti{2}1&=   1\anti{2}1 &\tau = 1\anti{2}1 &m = 0.\\        
    \end{array}
\end{equation*}

\begin{theorem} (Discrepancies)\label{thm:discrepancies}
    Suppose  $\ttt{T}$ is an expansion subtree of  $\tree(R)$.   
Let $Q = \rod{\inner(\ttt{T})}$ and $S = \rod{\leaves(\ttt{T})}$ be the rod sets corresponding to those node sets.
Then for $n>0$, $\discrepancySeq(n,R,S)$  is $\rodcount(n, Q)$, the net number of rods in $Q$ of length $n$.
\end{theorem}

\begin{proof}
Let $A$ be the subset of $\trains(R)$  for which $m= 0$
in Lemma~\ref{lemma:parseRtrains}, and let $B=\trains(R) \setminus A$.

For trains of length $n$, we have:
\begin{equation*}
    \tcount(n,R) =   \rodcount(n,A)+ \rodcount(n,B).
\end{equation*}
%where the sign of the train to its node determines the sign of the %train.

Each train in $B$ starts with the train $\sigma_1$, with
$\rod{\sigma_1} \in S$, and ends with an arbitrary element
of $\trains(  n - \text{length}(\sigma_1), R)$. Therefore
\begin{equation*}
\rodcount(n,B)= \sum_{k \in S} (\sign k)\tcount(n-k,R)
\end{equation*}
so
\begin{equation*}
\discrepancySeq(n,R,S) = \rodcount(n,A) .
\end{equation*}

Finally, $\rodcount(n,A)= \rodcount(n,Q)$
since the trains  $\tau \in A$ are precisely those that are nodes in $Q$.
\end{proof}

Now we can justify the claims we made in discussing the examples earlier in the section. Theorem~\ref{thm:discrepancies} implies Equation~\ref{eq:1anti2-6} since in the second tree in Figure~\ref{fig:qs-expansion} the  leaf rod set $\rod{\leaves(T)}$ consists of the length $6$ rod $ \rod{1\anti{2}1\anti{2}}$, the length $2$ rod/antirod pair $\rod{11}$/$\rod{\anti2}$ and 
the length $5$ rod/antirod pair  $\rod{1\anti{2}\anti{2}}$/$\rod{1\anti{2}11}$. 
Those pairs contribute $\pm 1$ when computing the net rod counts for the leaf nodes.
The internal nodes correspond to $\rod{1} ,\rod{1\anti{2}},$ and $\rod{1\anti{2}1}$, which are the discrepancies in \ref{eq:1anti2-6}. 

For Equations~\ref{eq:Fibidentity1} and \ref{eq:Fibidentity2} we use the third tree in Figure~\ref{fig:qs-expansion}, and note that the trains of odd length end at leaves and the trains of even length at internal nodes.

\section{Equivalence}

In this short section we introduce the first of the tools we need to construct our rod set algebra: formally defining how to cancel rod/antirod pairs as we did \emph{ad hoc} in the  the last section.

Different rod sets may have the same net train counts. For example, $\trains(\emptyrods)$ has just one element---the empty train with length $0$---while
$\trains([1, \anti{1}])$ is the infinite set of finite trains consisting of $1$'s and $\anti{1}$s. We can pair the elements of length $n$ in $\trains([1, \anti{1}])$ by switching the sign of the first rod, which guarantees that the trains in each pair have opposite signs, and shows that for any $n$ there are the same number of positive and negative trains, so 
$\tcount(n, [1, \anti{1}]) = \tcount(n, \emptyrods) $ for all $n$.

The next theorem shows that this kind of equality occurs only in the obvious way.

\begin{theorem}
    The following assertions about rod sets $R$ and $S$ are equivalent.

\begin{enumerate}
    \item $\rodcount(n,R) = \rodcount(n,S)$ for all $n > 0$.
    \item $\tcount(n,R) = \tcount(n,S)$ for all $n \ge 0$.
\item $R$ and $S$ can each be converted to the same rod set $T$ free of rod/antirod pairs by deleting some rod/antirod pairs $(k,\anti{k})$.

    \item $R$ can be converted to $S$ through additions and deletions of rod/antirod pairs $(k,\anti{k})$.

\end{enumerate}

\end{theorem}
\begin{proof}
The proof is an elementary consequence of the definitions; we leave it to the reader. 
\end{proof}

\begin{definition}
Rod sets $R$ and $S$ are \emph{equivalent}, written $R \equiv S$, when they satisfy any one (and hence all) of the conditions in the previous theorem. Each equivalence class contains a unique \emph{reduced} rod set $T(R)$ that contains no rod/antirod pairs. 
Rod set $R$ is \emph{positive} when $T(R)$ contains no antirods, is \emph{finite} when $T(R)$ contains finitely many rods, and is \emph{primitive} when the greatest common divisor of lengths of rods in $T(R)$ is $1$.
We write $\min R$ for the length of the shortest rod in $T(R)$.
When $R$ is finite we  write $\size{R}$ for the number of rods in $T(R)$ and $\max R$  for the length of the longest   rod in $T(R)$.
\end{definition}

For example, $R \equiv S$ when $R=\{1,\anti{1},2, 2, \anti{2},3,3\}$ and $S=\{2, 3,3,3,3,\anti{3},\anti{3},17,\anti{17}\}$. Then $T(R) = [2,3,3]$, $\size{R} = 3$, $\min R = 2$  and $\max R = 3$.

\begin{theorem}\label{thm:rodsetalgebra}
 If $R \equiv R'$ and $S \equiv S'$ then
 \begin{align*}
     R \cup S & \equiv R' \cup S' \\
          R \cap S & \equiv R' \cap S' \\
     RS &\equiv R'S' \\
     \anti{R} & \equiv    \anti{R'}
 \end{align*}
 Furthermore, if $R$ and $S$ are nonempty then

\begin{align}\label{eq:RbarS}
  R \cup \overline{R} &\equiv \emptyrods ,\notag \\
      \anti{(RS)} &\equiv \anti{R}S \equiv R\anti{S} \ne \anti{R}.\anti{S} ,\\
     \anti{R}.\anti{S} & \equiv RS. \notag 
 \end{align}
If $S \subseteq R$ then
    \begin{equation*}
        R \setminus S \equiv R \cup \anti{S}.
    \end{equation*} 
\end{theorem}

\begin{proof}
    The first seven assertions follow directly from the definitions. For the last one, write $R$ as the disjoint union of $R''$ and $S$. Then $R \setminus S=R''$, and $R \cup \anti{S}=R'' \cup S \cup \anti{S} \equiv R''$.
  \end{proof}
  
\begin{corollary}\label{cor:trainsR}

    The rod set $\rod{\trains(R)}$ is equivalent to the reduced rod set containing $\tcount(n,R)$ rods of length $n$ with appropriate signs.
\end{corollary}

\section{Rod Set Algebra}\label{sec: algebra}

In this section we prove that rod set expansion 
is well-defined up to rod set equivalence and is independent of which expansion tree we use. Then we develop algebraic ways to manipulate equivalence classes of rod sets. 
 
\begin{figure}
    \centering
    \includegraphics[valign=t]{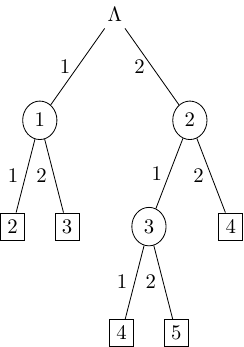}
  %  \phantom{aaaaa}
    \hspace{0.9in}
    \includegraphics[valign=t]{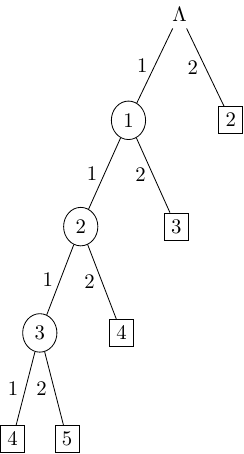}\\
       \includegraphics{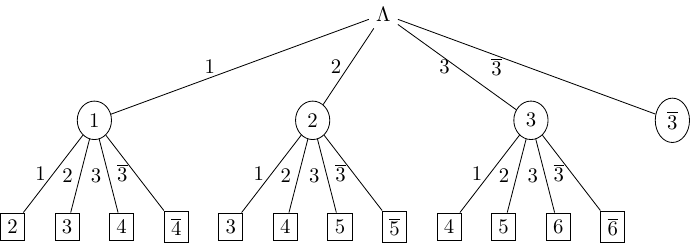}
    \caption{Three trees for the expansion
    $ [1,2]\trailright{[1,2,3]}[2,3,4,4,5]$.}
    \label{fig:threetrees}
\end{figure}

Figure~\ref{fig:threetrees}  shows three expansion trees for rod sets equivalent to $[1,2]$, with each node $\tau$ labeled with the rod $\rod{\tau}$. The three  leaf sets all satisfy $\rod{\leaves(\ttt{T})} \equiv [2,3,4,4,5]$, while the internal node sets satisfy
$\rod{\inner(\ttt{T})} \equiv [1,2,3]$. The figures suggest we can write
\begin{equation*}
[1,2]\trailright{[1,2,3]}[2,3,4,4,5]    
\end{equation*}
for any one of the three expansions
\begin{align*}
   \{1, 2\} &\trailright{\{1,2,21\}} \{11 ,12,211 ,212,22\},\\
    \{1, 2\} &\trailright{\{1,11,111\}} \{1111,1112,112,12,2\},\\
     \{1, 2, 3, \anti{3}\} &\trailright{\{1,2,3\}} 
     \{11 ,12,13,1\anti{3}, 
          21 ,22,23,2\anti{3}, 
        31 ,32,33,3\anti{3},\anti{3}                           
     \}.  
\end{align*}

 %valign

Since $R=[1,2]$, when we construct the first two trees in Figure~\ref{fig:threetrees} step-by-step we have to wait to expand a node with rod label $3$ until one appears in the tree. However, quantum mechanics suggests a way to avoid the wait: it predicts that systems (even the vacuum)
may spontaneously create particle/antiparticle pairs, which usually
mutually annihilate immediately. But sometimes one member of the pair
decays before annihilation happens. That metaphor motivates the
third expansion subtree in the figure, where we add the rod/antirod pair $[3,\anti{3}]$ to $[1,2]$ in order to expand rod $3$ without having to wait for it to show up. We end up with an equivalent leaf set.

Algebraically, the metaphor corresponds to adding and subtracting $\tcount(n-3, [1,2])$ in order to write the original recursion as if it were for the rod set $[1,2,3, \anti{3}]$:
\begin{equation*}
    \tcount(n,[1,2]) = \tcount(n,[1,2,3, \anti{3}])
    = \tcount(n-1) + \tcount(n-2) + \tcount(n-3) - \tcount(n-3),
\end{equation*}
 then using the recursion for $[1,2,3, \anti{3}]$ to replace the first three terms in the sum with four new terms each (two of these replacement terms also add to zero). After canceling equal terms what remains is the recursion for $[2,3,4,4,5]$.

If all the rods in $Q$ that we want to expand are already in $R$ then  the rod set $\rod{\leaves(\ttt{T})}$
for the expansion subtree $\ttt{T}$ of $\tree(R)$ with $\rod{\inner(\ttt{T})}=Q$ is 
$ (R \setminus Q) \cup QR $.
Theorem~\ref{thm:rodsetalgebra} then implies 
\begin{equation*}
  \rod{\leaves(\ttt{T})} = (R \setminus Q) \cup QR \equiv R \cup \anti{Q} \cup QR.
\end{equation*}

We will soon see that up to rod set equivalence
 this union is  the leaf rod set when expanding $R$ by $Q$, even when $Q$ is not a subset of $R$.

  \begin{definition}\label{def:rodexpansion}
  For rod sets $R$, $Q$ and $S$ we say $R$ expands to $S$ via $Q$, and write
  \begin{equation*}
      R \trailright{Q} S,
  \end{equation*}
  if there is a rod set $R' \equiv R$ and an expansion subtree $\ttt{T}$ of $\tree(R')$ for which $\rod{\inner(\ttt{T})} \equiv Q$ and $\rod{\leaves(\ttt{T})} \equiv S$ such that
     \begin{equation*}
    R' \trailright{\inner(\ttt{T})} \leaves(\ttt{T}) .
\end{equation*}
\end{definition}

 \begin{theorem} \label{thm:expansion} (Expansion) 
    For rod sets $R$ and $Q$ there is a rod set $S$ such that
    \begin{equation*}
        R \trailright{Q} S,
    \end{equation*}
    and this expansion is well-defined up to rod set equivalence: if $R' \equiv R$, $Q' \equiv Q$, and $\ttt{T}$ is an expansion subtree of $\tree(R')$ with $\rod{\inner(\ttt{T})} = Q'$, then
    \begin{equation*}
    \rod{\leaves(\ttt{T})} \equiv R \cup \anti{Q} \cup QR.
    \end{equation*}
\end{theorem}
\begin{proof}
We start by proving uniqueness.
Let $\ttt{T}$ be any expansion subtree of $R'\equiv R$ with internal nodes $Q' \equiv Q$. 
Every node in $\ttt{T}$ is either the root $\emptytrain$, a leaf, or an internal node in $Q'$. Furthermore, every node in $\ttt{T} \setminus \{\emptytrain \}$ is either a child of $\emptytrain$, hence a node in $R'$, or a child of an internal node, hence a node in $Q'R'$. 
Thus $R' \cup Q'R'$ is a list of all nodes in $\ttt{T} \setminus \{\emptytrain \}$. To list the leaves, we take out the internal nodes, 
\begin{align*}
   \rod{\leaves(\ttt{T})} &=\rod{(R' \cup Q'R')\setminus Q'} \\
    &\equiv (R' \cup Q'R') \cup \anti{Q'} \\
    &\equiv R \cup \anti{Q} \cup QR.
\end{align*}

    To show existence we need just one instance. Let $R' = R \cup \anti{Q} \cup Q \equiv R.$
Then
expanding $Q$ in $\tree(R')$ produces the expansion subtree 
$  \ttt{T} =  \{\emptytrain\}  \cup R'  \cup QR' $ for which
          $\inner(\ttt{T})= Q$.
          
% Then $  \ttt{T} =  \{\emptytrain\}  \cup R'  \cup QR' $
% is an expansion subtree of  
%           $\tree(R')$  for which
%           $\inner(\ttt{T})= Q$.
\end{proof}

Now we explore some simple consequences of this theorem  that enable 
useful algebraic operations on rod sets.

\begin{corollary}\label{cor:switchRQ} (Exchange) 
Rod set $R$ expands via $Q$ to $S$ if and only if 
  $\anti{Q}$ expands via $\anti{R}$ to $S$.
\end{corollary}
\begin{proof}
Since    $QR \equiv \anti{Q}.\anti{R}$,
\begin{equation*}
    R \cup \anti{Q} \cup QR = \anti{Q} \cup R  \cup \anti{Q}.\anti{R} .
\end{equation*}
\end{proof}

\begin{corollary}\label{cor:oddsignswap}
    (Odd Sign Swap) Suppose $R \trailright{Q} S$. Let $\oddss{R},\oddss{S}$, and $\oddss{Q}$ respectively represent $R,Q,S$ with the signs of odd length rods switched and those of even length rods kept the same. Then 
    \begin{equation*}
        \tcount(n, R^o) = (-1)^n \tcount(n, R) 
    \end{equation*}
    and
        \begin{equation*}
             \oddss{R} \trailright{\oddss{Q}} \oddss{S}.
    \end{equation*}
\end{corollary}
\begin{proof}
The first assertion follows from the fact that a train of odd (even) length contains an odd (even) number of rods of odd length. To prove the second, 
    we may assume $S=R \cup \anti{Q} \cup QR$. To show that $\oddss{S}=\oddss{R} \cup \anti{\oddss{Q}} \cup \oddss{Q}\oddss{R}$ we deal separately
 with each of the three rod sets whose union is $S$. Switching the signs of the odd length rods in $R$ and $\anti{Q}$ yields $\oddss{R}$ and $\anti{\oddss{Q}}$. The set $QR$ consists of rods of the form $\rod{q.r}$ with $q \in Q$ and $r \in R$. This concatenation has odd length when exactly one of $q$ and $r$ has odd length, so $\oddss{(QR)} = \oddss{Q}\oddss{R}$.
\end{proof}

The next theorem collects several important equivalent ways to characterize the expansion $\RQS$.

\begin{theorem}
\label{thm:augtrains}
Let $R$, $Q$ and $S$ be rod sets. Then the following assertions are equivalent
\begin{enumerate}
    \item 
\begin{equation*}
    \RQS .
\end{equation*}
\item For $n> 0$
 \begin{equation}\label{eq:DequalsQ}
      \discrepancySeq(n,R,S) = \rodcount(n,Q) .
\end{equation}
\item Trains built from rods in $R$ parse to trains built from $S$ and trains built from $S$ with a rod from $Q$ appended, so
\begin{equation}\label{eq:augtrains}
 \rod{\trains(R)} 
         \equiv \rod{\trains(S)} \cup  \rod{\trains(S).Q} .
\end{equation}
\item 
For $n > 0$
\begin{equation}\label{eq:countaugtrains}
\tcount(n,R) = \rodcount(n, \trains(R)) 
         = \rodcount(n,\trains(S)) + \rodcount(n,\trains(S).Q) .
\end{equation}
\end{enumerate}

\end{theorem} 

\begin{proof}
We use any convenient expansion subtree $\tree(R')$ for a rod set equivalent to $R$ and invoke Lemma~\ref{lemma:parseRtrains} and Theorem~\ref{thm:discrepancies}.
\end{proof}

As a corollary we see how to compose expansions. We write  $Q_{RS}$ when we need to name a rod set that expands $R$ to $S$. 

\begin{corollary}(Composition) 
    \label{cor:composition}
        Let $P,R,S$ be rod sets with $P \trailright {Q_{PR}} R$, \\
        $R \trailright {Q_{RS}} S$, and $P \trailright {Q_{PS}} S$. Then 
        \begin{equation*}
            Q_{PS} \equiv Q_{PR} \cup Q_{RS}  \cup Q_{PR}.Q_{RS}.
        \end{equation*}
\end{corollary}
\begin{proof}
     From the theorem we have,
\begin{align*}
    \rod{\trains(P)} &\equiv \rod{\trains(R) \cup \trains(R).Q_{PR}} \\
    \rod{\trains(R)} &\equiv \rod{\trains(S) \cup \trains(S).Q_{RS}}.
\end{align*} Substituting, we obtain 
    \begin{align*}
        \rod{\trains(P)} &\equiv \rod{\trains(S)\cup \trains(S).Q_{RS}} \\
        & \quad \quad  \cup \rod{(\trains(S) \cup \trains(S).Q_{RS}).Q_{PR}}.
    \end{align*}
    We distribute, view $Q_{RS}.Q_{PR}$ as one suffix, and the corollary follows.
\end{proof}

Theorem~\ref{thm:expansion} shows that there is a simple formula for $S$ in terms of $R$ and $Q$. The next theorem proves that given any two of $R,Q,S$, the third is well-defined. 

\begin{theorem}\label{thm:solveforQRS}
    Any two of $R,Q,S$ uniquely determine the third so that $R \trailright{Q} S$.
\end{theorem}

\begin{proof}
    Theorem~\ref{thm:expansion} tells us how to find $S$ given $R$ and $Q$. 

    Given $R$ and $S$  we can use Theorem~\ref{thm:augtrains} to find $Q$ by computing the discrepancies $\discrepancySeq(n, R,S)$. 

    Given $Q$ and $S$,  we again use use Theorem~\ref{thm:augtrains} to find the $\anti{R}$ for which $ \anti{Q}\trailright{\anti{R}} S$. Then Corollary~\ref{cor:switchRQ} implies $\RQS$. 
\end{proof}

It's interesting to examine edge cases.  When $Q$ is empty we have
            \begin{equation*}
        R \trailright{\emptyrods}  R.
    \end{equation*}

When $Q = \rod{\trains(R)}$
 every node in $\tree(R)$ has been expanded, there are no leaves, and
            \begin{equation*}
        R \trailright{\rod{\trains(R)}} \emptyrods .
    \end{equation*}   

When $R$ is empty, $ RQ = \emptyrods. Q =  \emptyrods$ is empty. If we add rod/antirod pairs $Q/\anti{Q}$ to  $[\ ]$ and expand everything in $Q$  then  no leaves other than $\anti{Q}$ remain, and we have 
            \begin{equation*}
        \emptyrods \trailright{Q}  \anti{Q} .
    \end{equation*}

\begin{figure}
    \centering
    \includegraphics[width=0.5\linewidth]{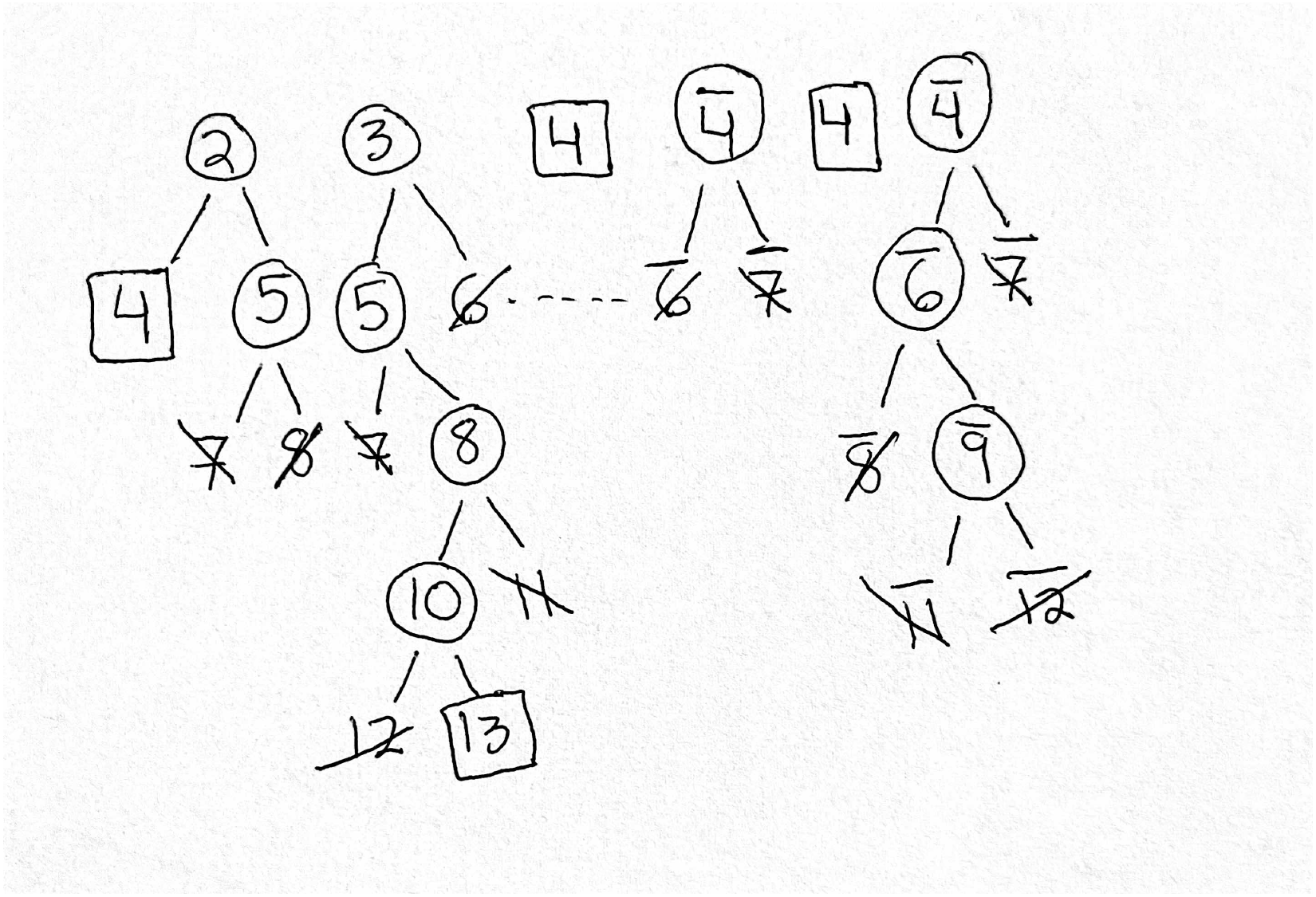}
    \caption{Hand calculation of $[2,3] \trailright{[2,3,\anti{4},\anti{4},5,5,\anti{6},8,\anti{9},10]} [4,4,4,13].$ Internal nodes, leaves, and canceled leaves are indicated respectively by circles, squares, and slashes. }
    \label{fig:handcalculation}
\end{figure}

In the proof of the preceding theorem we used discrepancies to find $Q$ given $R$ and $S$. In practice, we construct an expansion tree, and hence $Q$, incrementally by hand---a technique we invented well before we worked out the formal definitions of the concepts involved. Figure~\ref{fig:handcalculation} shows how we checked the surprising result of a computer search that told us that $[2,3]$ expands to $[4,4,4,13]$ via finite $Q$. 

We start by expanding both rods in $[2,3]$ since they do not appear in $S$. We draw a square around the resulting $4$ to indicate that it's a leaf. Since $S$ needs more $4$'s, and all the leaves we might expand are too big to have $4$ as a child, we add two $4/\anti{4}$ rod/antirod pairs at the root, then expand the $\anti{4}$'s and their descendants, canceling rod/antirod pairs as we go. Note that at length $8$ there are $2$ positive rods and one antirod. We can cancel either positive rod, leading to two non-isomorphic expansion trees. We then easily read the circled internal nodes $Q=[2,3,\anti{4},\anti{4},5,5,\anti{6},8,\anti{9},10]$ from the final picture.

Formally, the tree in the picture is not an expansion subtree of $\tree([2,3,4,\anti{4},4,\anti{4}])$ because many nodes are missing some children. For example, node 2 should have children 6, \anti{6}, 6 and \anti{6} in addition to 4 and 5. Since all missing nodes occur as rod/antirod pairs that are never expanded, we safely ignore them.

Below we formalize this algorithm. Follow along in
Figure~\ref{fig:handcalculation}.

\begin{algorithm}\label{alg:findQ}
    Let $R$ and $S$ be reduced rod sets with finite $S$. The following algorithm constructs
    a rod set $R' \equiv R$ and 
    a subtree $\ttt{T}$  of $\tree(R')$ for which
        $\rod{\leaves(\ttt{T})} \equiv S$
    and
    \begin{equation*}
            R \trailright{\rod{\inner(\ttt{T})}} S.
    \end{equation*}
 \end{algorithm}
    The algorithm constructs $R'$ and $\ttt{T}$ incrementally. At step $n$ we construct $R_n \equiv R$ and $\ttt{T}_n$,  a subtree of $\tree(R_n)$.  
    The leaves of $\ttt{T}_n$ are distributed among three sets: 
  
    \begin{itemize}
        \item $S_n = S \setminus (\leaves(\ttt{T}_n) \cap S$), the leaves we have not yet boxed in the figure. 
        \item $C_n$,  the canceled leaves, crossed out in the figure
        \item $E_n$, the remaining leaves, eligible for expansion.
    \end{itemize}
{\noindent
    Initialize:}
    \begin{itemize}
        \item  $R_0 = E_0=R$,
        \item $C_0=\emptyrods$,
          \item  $S_0 = S$, 
         \item $\ttt{T}_0$ is the expansion subtree of $\tree(R_0)$ with $\inner(\ttt{T}_0) = \emptyset$,
    \end{itemize}
{\noindent Loop on $n$:}

\phantom{foo}

The algorithm ends when one of the following conditions is true:

    \begin{itemize}
         \item $S_n = \emptyrods$. Then $R' = R_n$, $\ttt{T} = \ttt{T}_n$, and  $Q = \inner(\ttt{T})$ is finite.
        
        \item $S_n \subset R_n$. Then $R' = R_n$,   $\ttt{T}$ is the subtree of $\tree(R')$ built by  recursively expanding  the nodes in $E_n$ and their descendants, and $Q = \inner(\ttt{T})$ is infinite.
\end{itemize}

When neither of the above conditions is met,  let $s$ be one of the shortest rods in $S_n \setminus (S_n \cap R_n)$ and $e$ one of the shortest rods in $E_n$. Then update as follows: 

\begin{itemize}
\item If $e=s$ (both length and sign), then we have found another leaf to box:
\begin{itemize}
     \item $E_{n+1}=E_n \setminus \{e\}$, 
      \item $S_{n+1} = S_n \setminus \{s\}$,
     \item $R_{n+1} = R_n$, $C_{n+1} = C_n$, and $\ttt{T}_{n+1}=\ttt{T}_n$. 
 \end{itemize}
 \item If $\length(e) < \length(s)$, expand $e$:
 \begin{itemize}
     \item $\ttt{T}_{n+1}$ is the tree built by expanding $e$ to $e.R_n$ in $\ttt{T}_{n}$,
     \item $R_{n+1} = R_n $, $S_{n+1} = S_n$,
     \item $E_{n+1}$ is the reduced rod set corresponding to $(E_n \cup eR_n) \setminus \{e\}$,
     \item $C_{n+1}=C_n \cup ((E_n \cup eR_n \setminus \{e\}) \setminus E_{n+1})$ is the set formed by adding to $C_{n}$ the rod/antirod pairs that were canceled in reducing to form $E_{n+1}$.

 \end{itemize}
            \item If $\length(e) > \length(s)$ or if $e$ and $s$ have the same lengths, but different signs,   add the rod/antirod pair $s/\anti{s}$ to $R_n$, and expand $\anti{s}$: 
            \begin{itemize}
            \item $\ttt{T}_{n+1}$ is the tree built by 
             adding leaves $s$ and $\anti{s}$ as children of the root in $\ttt{T}_{n}$ and expanding $\anti{s} \rightarrow \anti{s}.R_n$,
                \item $R_{n+1} = R_n \cup \{s, \anti{s}\} $,
                \item $S_{n+1} = S_n \setminus \{s\}$,
                \item $E_{n+1}$ is the reduced rod set corresponding to $ E_n \cup \anti{s}R_n$,
                \item $C_{n+1}=C_n \cup ((E_n \cup \anti{s}R_n) \setminus E_{n+1})$ is the set formed by adding to $C_{n}$ the rod/antirod pairs that were canceled in reducing to form $E_{n+1}$.
            \end{itemize}
        \end{itemize}
\begin{proof}
    The algorithm terminates because one of the following will happen in a finite number of steps: $\size{S_n}$ decreases,  $\min E_n$ increases, or $\min S_n$ increases.

When  $S$ is infinite this algorithm may not terminate, but it does generate (the possibly infinite) $Q$ because the first two alternatives remove one of the shortest eligible rods. Thus for each $s \in S$ the third alternative must occur eventually.
\end{proof}

We have been careful to ensure that all our algebra works with infinite rod sets. It's clear that if $R$ expands to $S$ via $Q$ and both $R$ and $Q$ are finite, then so is $S$. The following theorems show some of what we can prove when we impose weaker finiteness conditions.

\begin{theorem}\label{thm:Qfinite}
    If $Q$ is finite and  $\RQS$, then for $n>\max Q$, $\trains(R)$ satisfies the recursion determined by $S$.
    Conversely, if for some $N$,  $\trains(R)$ satisfies the recursion determined by $S$ for $n>N$, then $Q$ is finite and $\max Q < N$.
    
    If $R$ is also finite, then $S$ is  finite and
        \begin{equation}\label{eq:addmax}
 \max S = \max R + \max Q.
    \end{equation}
\end{theorem}

\begin{proof}
    Equation~\ref{eq:DequalsQ} tells us that for $n>\max Q$, the discrepancy $\discrepancySeq(n,R,S)=0$.

    If $R$ is finite
then the longest rod in $R \cup \anti{Q} \cup QR$ comes from $QR$.
\end{proof}

For finite $R$ and $S$ the $Q$ in \RQS{} need not be finite.
  For example, when we use Algorithm~\ref{alg:findQ} 
to find the $Q$ that mediates the expansion of the Fibonacci rod set $[1,2]$ to the Padovan rod set $[2,3]$ we expand the $1$, add a $3/\anti{3}$ pair, then expand the $\anti{3}$ forever.
The result is
\begin{equation*}
Q = [1,2,3,4^2,5^3,6^5,7^8,8^{13}, \ldots].
\end{equation*}

The algorithm explains why the Fibonacci multiplicities appear in $Q$: expanding $\anti{3}$ forever means appending  $\anti{3}$ to every internal node in $\tree([1,2])$. The next theorem generalizes this observation.

\begin{theorem}\label{thm:RSfinite}
Suppose $R$ and $S$ are finite rod sets and $\RQS$. Let $m$ be the longest rod in $S\setminus R$. Then
 for  $n>m$,  
 $\rodcount(n,Q)$ satisfies the recursion determined by $R$. 
\end{theorem}

\begin{proof}
   Given $R$ and $S$, one way to construct a tree for the expansion is to add rod/antirod pairs to $R$ for the rods in $S$ that are not already in $R$, then expand everything that is not in $S$ forever. Formally, let $ R' = R \cup (S \setminus R) \cup \anti{(S \setminus R)} \equiv R$. Then $S \subseteq R'$, so $R'$ is a disjoint union $S \cup U$, with $\max U = m$.
   Let  \ttt{T} be the expansion subtree  of $\tree(R')$ for which $\inner(\ttt{T})$ is $U$ and all its descendants. 
      Then $S = \leaves(\ttt{T})$ and  $Q = \inner(\ttt{T})$.

$\inner(\ttt{T})$ is the disjoint union of the subtrees $\ttt{T}_k$ descending from the nodes $k \in U$. For each such $k$, the net number of nodes in $\ttt{T}_k$ of length $n$ is
$\tcount(n-k,R)$. For each $k$, the sequence $\tcount(n-k,R)$ satisfies the recursion for $R$ when $n > m \ge k$, so the same is true for their sum
    \begin{equation*}
        \rodcount(n,Q)=\sum_{k \in U}\tcount(n-k,R).
    \end{equation*}
 \end{proof}
 
\begin{corollary} \label{cor:maxR0inarow}
If $R$ and $S$ are finite rod sets with $\RQS$, then
 $Q$ is finite if and only if starting for some $N \ge \max S - \max R$ there are $\max R$ consecutive rod lengths $k$ for which $\rodcount(k,Q) = 0$.
\end{corollary}

\begin{proof}
    The forward direction is obvious. The converse follows from the fact that $\rodcount(n,Q)$ satisfies the recursion determined by $R$ for $n> \max S$ so all terms after $N$ are also zero, and thus $Q$ is finite.
\end{proof}

We close this section by using rod set algebra to study sequences that satisfy a finite recursion but don't start with our default initial conditions.

The Lucas numbers,
\begin{equation}\label{eq:lucasseq}
   \seq{L} = 2, 1,3, 4,7, 11, 18, \ldots  
\end{equation}
form a sequence that starts with $2,1$ and then
satisfies the Fibonacci recursion determined by the rod set $[1,2]$:
\begin{equation}\label{eq:lucas}
    \seq{L}(n) = \seq{L}(n-1) + \seq{L}(n-2)  .
\end{equation}

If we prepend $\seq{L}(0) = 1$ to the sequence in Equation~\ref{eq:lucasseq} to form 
\begin{equation}\label{eq:lucaswpre}
    1,\seq{L}=1,2,1,3,4,7,11,18, \ldots,
\end{equation} 
then Theorem~\ref{thm:recursionsarerodsets} implies there is a rod set $R$ whose net train counts are $1,\seq{L}$. But that rod set is uninformative: it does not clearly connect the Lucas numbers to the Fibonacci recursion determined by the rod set $[1,2]$. However, Theorem~\ref{thm:augtrains} says that the rod set $Q$ for which  
\begin{equation*}
    R \trailright{Q} [1,2]
\end{equation*}
is precisely the set of rods that specify when $1,\seq{L}$ differs from the Fibonacci sequence determined by $[1,2]$.

We start with $\seq{L}(1)=2$. If this train count satisfied the Fibonacci recursion, then it would be equal to $\seq{L}(0)+\seq{L}(-1)=1+0=1$, but it's one more than that, so the discrepancy at $n=1$ is $1$. We compute: 
\begin{equation*}
\begin{array}{rrrrrr}
    \rodcount(1,Q) &= \discrepancySeq(1,R,[1,2]) &= &2 - 1 - 0 &= &1,\\
    \rodcount(2,Q) &= \discrepancySeq(2,R,[1,2]) &= &1 - 2 -1 &= &-2.
\end{array}
\end{equation*}
For $n>2$,  $\seq{L}$ satisfies the Fibonacci recursion by definition, so we have $Q = [1, \anti{2}, \anti{2}]$.

Thus, also by Theorem~\ref{thm:augtrains}, the Lucas numbers satisfy the recursion
\begin{equation}\label{eq:lucasRec}
    \seq{L}(n)=\tcount(n,[1,2])+\tcount(n-1,[1,2])-2\tcount(n-2,[1,2]),
\end{equation}
which expresses each Lucas number in terms of Fibonacci numbers with our indexing and default initial conditions. 

Theorem~\ref{thm:augtrains} provides more. Equation~\ref{eq:countaugtrains} shows a nice way to view the trains counted by $\seq{L}$:
\begin{equation*}
 \rod{\trains(R)} 
         \equiv \rod{\trains([1,2])} \cup  \rod{\trains([1,2]).[1,\anti{2},\anti{2}]}. 
\end{equation*}
Table~\ref{table:lucas} illustrates the Lucas trains for the first few values of $n$ as either Fibonacci trains or as Fibonacci trains with an extra suffix from $Q$ .

Equation~\ref{eq:countaugtrains} constructs the extra trains using the rods in $Q$ as suffixes. Benjamin and Quinn\cite{benjamin} handle initial conditions a little differently.  They put the suffix first and call it a \emph{phase}. Their trains always start with a phase, so the trains in our model that don't have a suffix would correspond to using rods from $R$ as phases as well, and some of these rods might cancel with rods from $Q$.

 \begin{table}[h!]
     \centering
     \caption{Rod trains for the Lucas numbers}.
 
     \begin{tabular}{|r|r|r|r|}
\hline
\thead{Length} & \thead{Trains$([1,2])$} & \thead{Trains$([1,2]).[1, \anti{2}, \anti{2}]$} & \thead{Net Total} \\
     %\text{length} & $\trains([1,2])$ & $ \trains([1,2]).[1, \anti{2}, \anti{2}]$ &  net total  \\
     \hline
     1 & 1 &  $\emptytrain.1$ & 2 \\
     2 & 11,2                 & 1.1,$\emptytrain$.\anti{2},$\emptytrain$.\anti{2}                    & 1  \\
     3 & 111,12,21            &11.1,2.1,1.\anti{2},1.\anti{2}             &3   \\
     4 & 1111,112,121,211,22 
     & 111.1,12.1,21.1,11.\anti{2},2.\anti{2},11.\anti{2},2.\anti{2} & 4
      \\
     \hline
     \end{tabular}
     \label{table:lucas}
  \end{table}

The argument in this Lucas number example relies only  on Theorem~\ref{thm:augtrains}, not on the finiteness of $S=[1,2]$ nor on the fact that the discrepancies occur only for values $ n \le \max S$. But when these finiteness conditions are met the result shows how to use our rod model to define ``traditional" initial conditions that specify the first $n$ values of a sequence, then the rest of the sequence is defined by a recursion that looks back $n$ steps.

If 
$S$ is a finite rod set   and
    \begin{equation*}
        \seq{L} = \ell_1, \ell_2, \ldots
    \end{equation*}
    is a sequence starting at index $1$ for which
    \begin{equation*}
        \ell_n = \sum_{k \in S}(\sign k)\ell_{n-k} \quad \text{for } n> \max S ,
    \end{equation*}
    then  $\ell_1, \cdots, \ell_{\max S}$ are the \emph{traditional initial conditions} that together with $S$ determine $\seq{L}$.

\section{Generating functions}

In this section we see how our rod set model for counting trains offers a concrete combinatorial way to interpret the power series algebra that appears when we use generating functions to study recursions.
In the application sections that follow we will sometimes complement combinatorial rod set algebra arguments with references to  polynomials and power series.

Our algebraic operations on rod sets define interactions between
 two different formal power series: the rod count generating function and the recursion generating function. 

\begin{definition}\label{def:cpoly}
The \emph{rod count generating function} for rod set $R$ is the formal
power series 
\begin{equation*}
    \rodgeneratingfunction(x,R) = 
     \sum_{k \in R}(\sign k)x^k  = \sum_{n > 0}\rodcount(n, R)x^n.
\end{equation*}
Since there are no rods with length $0$, the constant term in this power series is $0$. When $R$ is finite, $\rodgeneratingfunction(x,R)$  is a polynomial.

  The \emph{recursion generating function} for the recursion defined by $R$ is the
  formal power series
  \begin{align*}
    \generatingfunction(x, R) &= \sum_{n=0}^\infty \tcount(n,R)x^n \\
    &= 1 + \tcount(1, R)x + \tcount(2, R)x^2 + \tcount(3, R)x^3 +
    \cdots \ .
\end{align*}
  \end{definition}

For example,
  \begin{equation*}
    \rodgeneratingfunction(x, [1,2]) = x + x^2,
\end{equation*}
and
  \begin{equation*}
\generatingfunction(x, [1,2]) = 1 + x + 2x^2 + 3x^3 + 5x^4 +     8x^5 + \cdots.
\end{equation*}

It follows directly from the definitions that 
        \begin{equation*}
         \rodgeneratingfunction(x,R) =   - \rodgeneratingfunction(x,\anti{R}),
    \end{equation*} 

    and that
    \begin{equation*}
                \generatingfunction(x,R) =  1+ \rodgeneratingfunction(x,\rod{\trains(R)}) .
    \end{equation*}  

    If $R \equiv R'$ then
    \begin{equation*}
        \rodgeneratingfunction(x,R) =  \rodgeneratingfunction(x,R') 
        \text{  and  } 
                \generatingfunction(x,R) =  \generatingfunction(x,R'). 
    \end{equation*}  

The next lemma establishes basic properties of the rod count generating function. 

\begin{lemma}
    For rod sets $Q,R$,
       \begin{equation}\label{eq:unionC}
        \rodgeneratingfunction(x,Q \cup R) =
          \rodgeneratingfunction(x,Q) +  \rodgeneratingfunction(x,R)
    \end{equation}  
    and
    \begin{equation}\label{eq:productC}
        \rodgeneratingfunction(x,QR) =
          \rodgeneratingfunction(x,Q)   \rodgeneratingfunction(x,R).
    \end{equation}
\end{lemma}
\begin{proof}
The first assertion is obvious. To prove the second, note that
    the rod set $QR$ consists of the rods $\rod{qr}$ with $q \in Q$ and $r \in R$, which correspond to the terms in the sum that defines
    $\rodgeneratingfunction(x,QR)$.
\end{proof}

The algebra for rod count and recursion generating functions works smoothly because we defined the rod count generating function
$\rodgeneratingfunction(n,R)$ only for $n>0$. We add or subtract this quantity from $1$ as needed, as illustrated in the next theorem. 

\begin{theorem}\label{thm:productformula}
    Rod set $R$ expands to $S$ via $Q$ if and only if
    \begin{equation}\label{eq:CS-CQCR}
         1 - \rodgeneratingfunction(x,S) = (1-\rodgeneratingfunction(x,R))(1+\rodgeneratingfunction(x,Q)).
    \end{equation}
\end{theorem}

\begin{proof}
If $R$ expands to $S$ via $Q$ then
    Theorem~\ref{thm:expansion} implies $S \equiv R \cup \anti{Q} \cup QR$ so
    \begin{align*}
        1 - \rodgeneratingfunction(x,S) &=
        1 - (\rodgeneratingfunction(x,R) -
        \rodgeneratingfunction(x,Q)
        + \rodgeneratingfunction(x,QR))\\
        &=       1 - (\rodgeneratingfunction(x,R) -
        \rodgeneratingfunction(x,Q)
        + \rodgeneratingfunction(x,Q)\rodgeneratingfunction(x,R))\\
        &= (1-\rodgeneratingfunction(x,R))(1+\rodgeneratingfunction(x,Q)).
    \end{align*}
    Conversely, this equation implies $S \equiv R \cup \anti{Q} \cup QR$, so
    $R$ expands to $S$ via $Q$.
\end{proof}

Equation~\ref{eq:CS-CQCR} shows how useful it is to define  rod generating functions without specifying $\rodgeneratingfunction(0,R)$: we need both of
$1 \pm \rodgeneratingfunction(x,R)$. The
  commutativity of polynomial multiplication does not smoothly translate to the rod context, where $R$ and $Q$ play different roles. For example, the polynomial equation 
\begin{equation*}
    (1-x-x^3-x^4)=(1-x-x^2)(1+x^2), 
\end{equation*}
corresponds both to the 
Fibonacci expansion $[1,2]\trailright{[2]}[1,3,4]$ and to its exchange, $[\anti{2}]\trailright{[\anti{1}, \anti{2}]}[1,3,4]$. For the first expansion, we derive the right side of the equation from $\rodgeneratingfunction(x,R)=x+x^2$ and $\rodgeneratingfunction(x,Q)=x^2$, and we derive the  second  from $\rodgeneratingfunction(x,R)=-x^2$ and $\rodgeneratingfunction(x,Q)=-x -x^2$.

The infinite Fibonacci expansion $[1,2]\trailright{[\text{Evens}]}[\text{Odds}]$ corresponds to the power series identity 
\begin{equation*}
    (1-x-x^3-x^5-\cdots)=(1-x-x^2)(1+x^2+x^4+\cdots).
\end{equation*}

The Exchange Corollary~\ref{cor:switchRQ} is an immediate consequence of Theorem~\ref{thm:productformula}. Odd sign swap (Corollary~\ref{cor:oddsignswap}) follows when we replace $x$ with $-x$ in Equation~\ref{eq:CS-CQCR}. 

This theorem shows that the Expansion Theorem~\ref{thm:expansion} is essentially equivalent to the standard algorithm for multiplying formal power series. The net count of the rods of length $n$ in $R \cup \anti{Q} \cup QR$ encodes the calculation of the coefficient of $x^n$ in the product power series. Composition (Corollary~\ref{cor:composition}) is just the distributive property.
Algorithm~\ref{alg:findQ}, illustrated in Figure~\ref{fig:handcalculation}, calculates the quotient of two power series.

While working on this paper we often found that we made fewer arithmetic errors when we calculated expansions using trees  
rather than by multiplying and dividing polynomials or power series using the standard algorithms. However, we do not recommend teaching our methods to children.

\begin{corollary}\label{cor:z[x]}
    If $R$ is a finite rod set then the rod count generating functions for the rod sets $S$ that are expansions of $R$ by a finite $Q$ are precisely the polynomials  with   constant term $1$ in the ideal generated by
    $1 - \rodgeneratingfunction(x,R)$ in the polynomial ring $\mathbb{Z}[x]$. 
\end{corollary}
\begin{proof}
Every polynomial with constant term $1$ is $1-\rodgeneratingfunction(x,R)$ for a finite rod set $R$, and conversely.
    If polynomial $s(x)$  has constant term $1$ and satisfies
    \begin{equation*}
    s(x) = q(x)(1- \rodgeneratingfunction(x,R))
    \end{equation*}
    then $q(x)$ has constant term $1$ and
    $\RQS$ for the rod sets $S$ and $Q$ corresponding to the polynomials $s(x)$ and $q(x)$.
\end{proof}

The next corollary is a reprise of Theorem~\ref{thm:solveforQRS}.

\begin{corollary}\label{cor:RSimpliesQ}
    Any two of $R,Q,S$ uniquely determine the third so that $R \trailright{Q} S$.
\end{corollary}
\begin{proof}
Using Equation~\ref{eq:CS-CQCR} we can solve for any one of the rod count generating functions in terms of the other two.
\end{proof}

\begin{theorem}\label{thm:generatingfunction}
  For rod set $R$
  \begin{equation}\label{eq:generatingfunction}
    \generatingfunction(x,R) =
    \frac{1}{ 1 - \rodgeneratingfunction(x, R)},
\end{equation}
and $R$ expands to $S$ via $Q$ if and only if
  \begin{equation}\label{eq:ctimesc}
    \generatingfunction(x,R)
 = (1+\rodgeneratingfunction(x,Q))\generatingfunction(x,S).
  \end{equation}
\end{theorem}
\begin{proof}
Since $R$ expands to $S = \emptyrods$ via $Q = \rod{\trains(R)}$, and $\rodgeneratingfunction(x,\emptyrods) = 0$,
 Equation~\ref{eq:CS-CQCR} implies
\begin{align*}
 1 &= 1 - \rodgeneratingfunction(x,\emptyrods) \\
&= (1 -  \rodgeneratingfunction(x,R))(1+ \rodgeneratingfunction(x,\rod{\trains(R)})) \\
&= (1 -  \rodgeneratingfunction(x,R))\generatingfunction(x,R).
\end{align*}
The second assertion follows from Theorem~\ref{thm:productformula}.
\end{proof}

For example, the recursion generating function for the Fibonaccis is 
\begin{equation*}
    \generatingfunction(x,[1,2])=\dfrac{1}{1-x-x^2}.
\end{equation*}
For the expansion $[1,2]\trailright{[2]}[1,3,4]$, we have 
\begin{equation*}
    \generatingfunction(x,[1,2])=\dfrac{1}{1-x-x^2}=\dfrac{1+x^2}{1-x-x^3-x^4}=(1+x^2)\generatingfunction(x,[1,3,4]).
\end{equation*}

The formal expansion of the geometric series in
Equation~\ref{eq:generatingfunction} is
  \begin{align*}
     \frac{1}{ 1 - \rodgeneratingfunction(x, R)} &=
     1 + \rodgeneratingfunction(x, R) +  \rodgeneratingfunction(x, R)^2 + \cdots\\
     &=    1 + \rodgeneratingfunction(x, R) +  \rodgeneratingfunction(x, R.R) + \cdots. 
\end{align*}
which nicely mirrors the decomposition of $\trains(R)$ in Equation~\ref{eq:geometrictrains}.

The next corollary is a reprise of Theorem~\ref{thm:RSfinite}.

\begin{corollary}
    Suppose $R$ and $S$ are finite rod sets and $\RQS$. Then 
 for \\ $n>\max S$,  
 $\rodcount(n,Q)$ satisfies the recursion determined by $R$. 
\end{corollary}
\begin{proof}
    Equations~\ref{eq:CS-CQCR}  and \ref{eq:generatingfunction}
    imply that 
    \begin{equation*}
        1+\rodgeneratingfunction(x,Q) = 
        (1-\rodgeneratingfunction(x,S))\generatingfunction(x,R). 
    \end{equation*}
    The right side of the equation is a finite sum of shifts of the net train count generating function for $R$. For $n>\max S$, there are no more new shifts, and the coefficients of $\rodcount(n,Q)$ satisfy the recursion for $R$.
\end{proof}

The following theorem is a reprise of Theorem~\ref{thm:recursionsarerodsets}.  

\begin{theorem}
    Every formal power series with constant term $1$ is the recursion
    generating function  
    $\generatingfunction(x,R)$ for a rod set $R$ unique up to rod/antirod equivalence.
\end{theorem}
\begin{proof}
  Let $\tcount(x) =  1 + \sum_{k=1}^\infty a_kx^k$.
  We calculate the formal power series inverse 
\begin{equation*}
\frac{1}{\tcount(x)} =  1 + \sum_{k=1}^\infty b_kx^k
\end{equation*}
Then we let $R$ be the rod set with $b_k$ rods (or antirods)
of length $k$ .
\end{proof}

We end this section with a well-known theorem on rational generating functions~\cite[p 202]{stanley}:

\begin{theorem}\label{thm:rationalGF}
The net train counts for a rod set $R$ eventually satisfy the recursion for a finite rod set $S$ if and only if $\generatingfunction(x,R)$ is a quotient of polynomials.
\end{theorem}
\begin{proof}
Suppose $\rodcount(n,R)$ eventually satisfies the recursion determined by a finite rod set $S$.
Let  $Q$ be the rod set for which
\begin{equation*}\label{eq:initialsgeneratingfunction}
    R \trailright{Q} S.
\end{equation*}
Since $Q$ specifies the discrepancies in this expansion, Theorem~\ref{thm:Qfinite} implies it is finite,
so
\begin{equation}\label{eq:rationalfunction}
    \generatingfunction(x,R) =
    \frac{1 + \rodgeneratingfunction(x,Q) }{1 - \rodgeneratingfunction(x,S)}
\end{equation}
 is a quotient of polynomials.

Conversely, if the numerator and denominator in Equation~\ref{eq:rationalfunction} are polynomials, then $Q$ and $S$ are finite and the net train counts for $R$ satisfy the recursion for $S$ when $n > \max Q$.
\end{proof}

Recall the expansion~\ref{eq:lucaswpre}
\begin{equation*}
    R \trailright{[1, \anti{2}, \anti{2}]} [1,2],
\end{equation*}
where $R$ is the rod set for the Lucas numbers, $1,\seq{L}$. By Equation~\ref{eq:rationalfunction}, the generating function for this sequence is
\begin{equation*}
    \frac{1+x-2x^2}{1-x-x^2} = 1 + 2x + x^2 + 3x^3 + 4x^4 + 7x^5 + \cdots .
\end{equation*}

To find the generating function for the Lucas numbers with traditional initial conditions \ref{eq:lucasseq},  we can subtract $1$ and divide by $x$:
\begin{equation*}
    \frac{2-x}{1-x-x^2} = 2 + x + 3x^2 + 4x^3 + 7x^4 + \cdots .
\end{equation*}

\section{Duality and its Applications}\label{sec:duality}

When  $R$ expands to $S$ via $Q$, the \emph{dual} of $Q$ is the rod set that expands $S$ to $R$. Duality is especially useful for exploring relationships between combinatorial identities and compositions. We end this section with Fibonacci and Padovan examples.

\begin{definition}
    The \emph{dual} of rod set $Q$ is
    \begin{equation*}
        Q^* = \rod{\trains(\anti{Q})}.
    \end{equation*}
\end{definition}

\begin{theorem}\label{thm:duality}
If
\begin{equation*}
    R \trailright{Q} S,
\end{equation*}
then
\begin{equation*}
    S \trailright{Q^*} R .
\end{equation*}

We call this the \emph{dual expansion} of $R \trailright{Q} S$.
\end{theorem}

\begin{proof}
 Figure~\ref{fig:treeduality} shows an expansion of $\tree(R \cup \anti{Q} \cup QR)$. Each node in the picture is actually a set of siblings corresponding to the rods in the rod set label. The powers of $\anti{Q}$ are the internal nodes.

\begin{figure}[h]
    \centering
    \includegraphics[width = 0.7\textwidth]{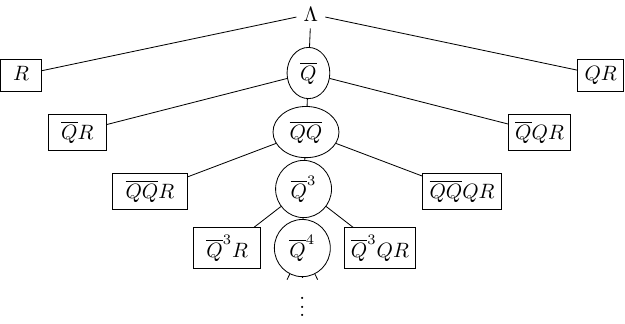}
    \caption{Expanding $\tree(R \cup \anti{Q} \cup QR)$ via $\trains \anti{Q}$  }
    \label{fig:treeduality}
\end{figure}
The union  of the internal nodes is $\trains(\anti{Q})$. The leaves on the left are $\anti{Q}^kR$ starting at $k=0$, while the leaves on the right are 
$\anti{Q}^{k-1}QR$, starting at $k=1$. Equation~\ref{eq:RbarS} shows that all but the leaves in $R$ cancel.
\end{proof}

\begin{corollary}
      For rod set $Q$, 
    \begin{equation*}
        Q^{**} \equiv Q.
    \end{equation*}  
\end{corollary}
\begin{proof}
    If $\RQS$, then $S \trailright{Q^*} R$ so
    $R \trailright{Q^{**} } S$. Then Theorem~\ref{thm:expansion} implies $Q^{**} \equiv Q$.
\end{proof}

 It is easy to write the rod generating function for the dual $Q^*$ in terms of the rod generating function for $Q$.

\begin{theorem}\label{cor:1overQ}
For rod set $Q$, we have 
\begin{equation*}
    1+\rodgeneratingfunction(x, Q^*) = 
    \frac{1}{1+\rodgeneratingfunction(x,Q)}.
\end{equation*}
\end{theorem}
\begin{proof} 
The dual $Q^{*}$ counts the trains of $\anti{Q}$ so      \begin{equation*}
    1+\rodgeneratingfunction(x, Q^*) = 
    \generatingfunction(x, \anti{Q}) =
    \frac{1}{1-\rodgeneratingfunction(x,\anti{Q})}=\frac{1}{1+\rodgeneratingfunction(x,Q)}.
\end{equation*}
\end{proof}

 When $\RQS$, Theorem~\ref{thm:expansion} shows how to compute $S$ given $R$ and $Q$, and Theorem~\ref{thm:solveforQRS} tells us that  any two of $R,Q,S$ determines the third.
 With duality, it is straightforward, albeit a bit messy, to find  explicit rod set formulas for $Q$ in terms of $R$
  and $S$ and for $R$ in terms of $Q$ and $S$.   

 Duality is particularly useful when expanding by a
 rod set $Q$ containing just a single rod $q$. Then $\trains(\anti{Q})$ contains rods with lengths that are all the multiples of the length of $q$. Duality then reveals an interesting connection between identities involving sums and those involving compositions.

 For example, recall the Fibonacci identities in
Equations~\ref{eq:Fibidentity1} and \ref{eq:Fibidentity2}, 
\begin{align*}
     1 + 2 + 5 + 13 +  34 +  \ \cdots   \   +\tcount(2n)    &= \tcount(2n+1) ,\\
    1 + 3 + 8 + 21 + \cdots + \tcount(2n-1) &= \tcount(2n) -1 ,
\end{align*}
which we established 
using the third expansion tree in Figure~\ref{fig:qs-expansion}. With Theorem~\ref{thm:augtrains}, we can derive these identities from the infinite expansion 
% For example, recall that we established the Fibonacci identities in
% Equations~\ref{eq:Fibidentity1} and \ref{eq:Fibidentity2}
% using the third expansion tree in Figure~\ref{fig:qs-expansion} together with Theorem~\ref{thm:augtrains}. We started with the infinite expansion
\begin{equation} \label{eq:inffib}
    [1,2]\trailright{\rod{\trains([2])}}[1,3,5,7,\ldots].
\end{equation}

Now we use Theorem~\ref{thm:duality} to dualize Equation~\ref{eq:inffib}:
\begin{equation} \label{eq:135...toFib}
    [1,3,5,7,\ldots]\trailright{[\anti{2}]}[1,2].
\end{equation}

By Theorem~\ref{thm:Qfinite}, the  train counts for the rod set $[1,3,5,7,\ldots]$ 
satisfy the Fibonacci recursion when $n>2$. The train counts for the rod set of odd numbers are the 
number of compositions, $\seq{O}(n)$, of $n$ into odd parts.  Since   $\seq{O}(1) = \seq{O}(2) = 1$, the sequence $\seq{O}(n)$ is the Fibonacci sequence with these initial conditions, hence $\seq{O}(n) = \tcount(n-1,[1,2])$. For example, there are $\tcount(4,[1,2])=5$ rod trains of length $5$ built from odd length rods: $5,311, 131, 113, 11111$.  

We can create similar relationships between
identities and compositions by starting with a rod set $R$ containing $q$ and expanding along all the edges 
in $\tree(R)$  with label $q$ to create an expansion with internal nodes
$Q=\rod{\trains([q])}$. The leaves attached to these internal nodes form a union of arithmetic progressions with difference $q$. The dual expansion has $Q^*=[q]$, with $\size{Q^*}=1$, and it reveals an identity about compositions.

Here is a Padovan example, expanding
$[2,3]$  by  $\rod{\trains([2])}$:
\begin{equation*} 
    [2,3]\trailright{\rod{\trains([2])}}[3,5,7,9,\ldots],
\end{equation*}
with dual expansion
\begin{equation*}
    [3,5,7,9,\ldots]\trailright{[\anti{2}]} [2,3].
\end{equation*}

Letting $P(n)=\tcount(n,[2,3])$, the original expansion leads to  the identities
\begin{equation*}
    \begin{array}{ll}
 P(0)+P(2)+P(4)+\cdots + P(2n)&=P(2n+3)   \\
 P(1)+P(3)+P(5)+\cdots + P(2n+1)&=P(2n+4) -1.
    \end{array}
\end{equation*}
From the dual expansion we see that the number of compositions of $n \ge 3$ into odd parts greater than or equal to $3$ satisfies the Padovan recursion. Examining the initial conditions tells us that the number of such compositions is  $P(n-3)$. 

 If instead we expand the Padovans by $\rod{\trains([3])}$, we discover
 \begin{equation*}
    \begin{array}{lll}
 P(0)+P(3)+P(6)+\cdots + P(3n)&=P(3n+2)   \\
 P(1)+P(4)+P(7)+\cdots + P(3n+1)&=P(3n+3) -1 \\
 P(2)+P(5)+P(8)+\cdots + P(3n+2)&=P(3n+4) .
    \end{array}
\end{equation*}

Then the dual expansion shows that the number of compositions of $n \ge 4$ into parts congruent to $2 \pmod{3}$ satisfies the Padovan recursion. Examining the initial conditions tells us that the number of such compositions is  $P(n-2)$. The reader is invited to find more such identities. 

The expansion $[\anti{1},\anti{2},\anti{3}\ldots]\trailright{[\anti{1}]} \emptyrods$ tells us that for $n > 1$, the number of compositions of $n$ into an odd number of parts is equal to the number of compositions of $n$ into an even number of parts. The dual expansion is less interesting in this case. 

The expansion 
\begin{equation}\label{eq:expand123dots}
    [1,2,3,\ldots]\trailright{[\anti{1}]} [1,1]
\end{equation}
  connects two well-known theorems. The forward expansion tells us that the total number of compositions of $n$ is $2^{n-1}$, while the dual expansion gives the identity  
  \begin{equation}\label{eq:sumpowersof2}
  1+2+2^2+\cdots +2^{n-1}=2^n-1.        
  \end{equation}

We end this section with a more complicated application using an expansion of two copies of an arithmetic sequence, 
\begin{equation*}
  S =   [2,2,5,5,8,8,11,11,\ldots]\trailright{[\anti{3}]}[2,2,3],
\end{equation*}
with dual expansion
\begin{equation*}
    R = [2,2,3]\trailright{\rod{\trains([3])}} [2,2,5,5,8,8,11,11,\ldots].
\end{equation*}

Since $R = [2,2,3]$ is an expansion of the Fibonacci rod set $[1,2]$  by the singleton rod set $[1]$, Corollary~\ref{cor:composition} implies
$[1,2] \trailright{}S$ has discrepancies $[3,6,9,\ldots] \cup [1,4,7,\ldots]$. The discrepancies are the rods with lengths $n \equiv 0,1 \pmod{3}$, which leads to the

following identities for the Fibonacci sequence $\tcount(n)=\tcount(n,[1,2])$:
\begin{equation*}
    \begin{array}{ll}
    \tcount(0)+\tcount(3)+\tcount(6)+\cdots +\tcount(3n)&=\dfrac{\tcount(3n+2)}{2}, \\[8pt]  
\tcount(1)+\tcount(4)+\tcount(7)+\cdots +\tcount(3n-2)&=\dfrac{\tcount(3n)-1}{2}, \\[8pt]   
\tcount(2)+\tcount(5)+\tcount(8)+\cdots +\tcount(3n-1)&=\dfrac{\tcount(3n+1)-1}{2}.
    \end{array}
\end{equation*}

Thus we can also use rod set algebra to derive identities that require fractions!

\section{Rod Sets and Binomial Coefficients}

It's well-known that the Fibonacci numbers are sums of diagonals in Pascal's triangle\footnote{a.k.a. Yang Hui's, Tartaglia's, \emph{et.al}'s triangle}.
More generally, when $R$ consists of two positive rods, the net train counts $\tcount(n,R)$ also appear as sums of diagonals, which leads to some nice binomial identities.  
Even more generally,  
when $R$ is positive and finite we can (but won't) systematically compute $\tcount(n,R)$ as a sum of multinomial coefficients on hyperplanes in high dimensional analogues of the triangle.

   \begin{definition}
        Let $R$ be a rod set and $\tau \in \trains(R)$. The \emph{underlying partition} of $\tau$ is the multiset of rods in $\tau$. 
    \end{definition}

    For example, if $R=[3,4]$ and $\tau=34443$ then the underlying partition of $\tau$ is $\{3,3,4,4,4\}$.

     \begin{theorem}\label{thm:binomtcount}
        Let $R=[a,b]$ be a positive rod set. Then 
        \begin{equation}\label{eq:binomsumgen}
            \tcount(n,R)= \sum_{\substack{ia+jb=n \\ i,j \ge 0}}\binom{i+j}{j}.
        \end{equation}
    \end{theorem}

    \begin{proof}
       For each pair of integers $i,j \ge 0$ with $ia+jb=n$, let $U$ be the multiset containing $i$ rods of length $a$ and $j$ rods of length $b$. A length $n$ train in $\trains(R)$ with underlying partition $U$ consists of $i+j$ rods for some pair $i,j$. We count the number of these trains by starting with $i+j$ places and then choosing in $\binom{i+j}{j}$ ways the places for a rod of length $b$. We sum over the possible partitions to calculate $\tcount(n,R)$.
    \end{proof}

   If $R$ contains more than two rods, we can use the same logic to write $\tcount(n,R)$ as a sum of multinomial coefficients; the linear constraint will just be an equation with more variables. When $R$ contains antirods we can adjust the signs in the summation. 

    If $R=[a,b]$ is primitive and positive, then the Chinese Remainder Theorem implies that integer solutions $(i,j)$ to $ia+jb=n$, the condition on the sum in Equation~\ref{eq:binomsumgen}, are congruent$\mod{a \times b}$. We care only about solutions with nonnegative $i$ and $j$. 
    
    To visualize the Pascal triangle diagonals, we order the terms in the sum  by increasing $j$, starting with the fewest possible rods of length $b$. At each step we replace  $b$  rods of length $a$ with $a$ rods of length $b$ to build the next partition, ending with the partition with the most rods of length $b$.

\begin{figure}
    \centering
    \includegraphics[width=\linewidth]{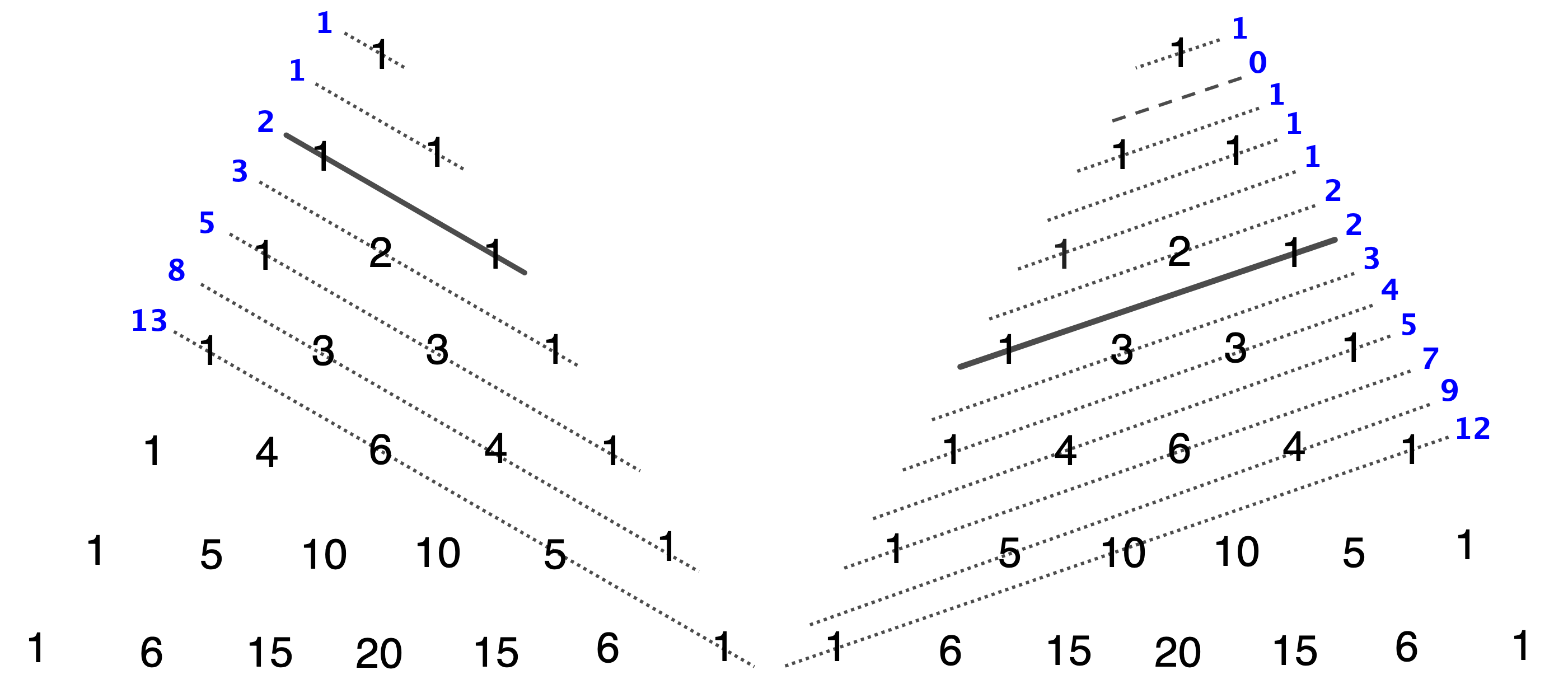}
    \caption{Train counts for the Fibonaccis, $[1,2]$, and Padovans, $[2,3]$, as diagonal sums in Pascal's triangle   }
    \label{fig:triangles}
\end{figure}

    For example, if $R=[3,5]$ and $n=70$, we start with the underlying partition $\{3^{20},5^2\}$ and then move to the underlying partition $\{3^{15},5^5\}$, and continue until we have the sum
    \begin{equation}\label{eq:binomsum}
        \tcount(70,[3,5])=\binom{22}{2}+\binom{20}{5}+\binom{18}{8}+\binom{16}{11}+\binom{14}{14}=63862.
    \end{equation}
  As we move from term to term in the sum, the number of rods in the underlying partition decreases by $b-a=2$, and the number of rods of length $b$ increases by $a=3$. The constant increases and decreases lead to the terms in the sum appearing on parallel diagonals in the triangle.
  Figure~\ref{fig:triangles} shows the diagonals for the Fibonacci and Padovan rod sets.
  When 
 $R=[a,b]$ is primitive, the first diagonal with two terms is
 \begin{equation*}
     \tcount(a \times b, [a,b]) =\binom{a}{0} + \binom{b}{b} =\binom{a}{a} + \binom{b}{0} = 2.
 \end{equation*}
 The two ways of writing this sum correspond to the dark diagonals in Figure~\ref{fig:triangles}. The first sum counts rods of length $b$ and the second rods of length $a$.  In practice, it's hard to draw the diagonals for large $a$. 
 Hopkins \cite[Appendix A]{hopkins} shows similar diagrams.

Theorem~\ref{thm:binomtcount} leads to some nice binomial identities. The train counts for $R=[1,1]$ are powers of two. The corresponding ``diagonals" are the rows of the triangle, providing another proof of the familiar identity 
     \begin{equation*}
     \tcount(n, [1,1]) =    \sum_{i=0}^n \binom{n}{i}=2^n.
     \end{equation*} 

When we allow antirods in $[a,b]$ and adjust the signs of the binomial coefficients appropriately we see more nice identities.

For $R=[1,\anti{1}]$   the net train counts are all zero for $n \ge 1$, hence
      \begin{equation*}
       \tcount(n, [1,\anti{1}]) =  \sum_{i=0}^n (-1)^i\binom{n}{i}=0,
     \end{equation*} 
     where $i$ represents the number of antirods in the underlying partition. 

When $R=[\anti{1},\anti{2}]$ the total train counts are the Fibonacci numbers, but the net train counts have period $3$.
Checking the signs carefully leads to the 
identity
    \begin{equation}\label{eq:period3}
          \tcount(n, [\anti{1},\anti{2}]) =   \sum_{0\le i \le n/2}(-1)^{n-i}\binom{n-i}{i}= 
        \begin{cases}
              1  & n \equiv 0 \pmod{3}\\
               -1 & n \equiv 1 \pmod{3} \\
                0 & n \equiv 2 \pmod{3}.
       \end{cases}
    \end{equation}
    Here $i$ is the number of rods of length $2$ in the underlying partition and $n-i$ is the total number of rods. Since both rods are antirods, the sign depends on the parity of the number of rods.

\section{Periodicity}\label{sec:periodicity}

At the end of Section~\ref{sec:duality} we looked at expansions where $\size{Q}=1$; here we'll explore expansions where $\size{S}=1$. These seemingly similar conditions   are conceptually different: the latter leads to periodicity.
For example, the expansion
\begin{equation}\label{eq:period3expansion}
    [\anti{1}, \anti{2}] \trailright{[\anti{1}]} [3]
\end{equation}
implies that $\tcount(n,[\anti{1}, \anti{2}])=\tcount(n-3,[\anti{1}, \anti{2}])$ for $n \ge 3$, which
explains the periodicity in the train counts
\begin{equation}
    \tcount(n,[\anti{1}, \anti{2}]) = 1, -1, 0, 1, -1, 0, \ldots 
\end{equation}
that we found at the end of the last section.

This example generalizes easily.

\begin{definition}
    The net train counts for rod set $R$ are \emph{eventually periodic} with \emph{period} $p> 0$ when there is an $N\ge 0$ such that
    \begin{equation*}
        \tcount(n+p, R) =   \tcount(n,R) \text{ for all } n \ge N.
    \end{equation*}

When $N=0$ we say the sequence is \emph{periodic}.

    Any positive multiple of a period is a period. Sometimes we will use \emph{the period} for the least period.
\end{definition}

The following theorem characterizes rod sets with eventually periodic net train counts.

\begin{theorem}(Periodicity)\label{thm:periodicity}
    The net train counts for rod set $R$ are eventually periodic with period $p$
   if and only if $R$ expands to $[p]$ via a finite $Q$. 
   The net train counts $\tcount(n,R)$ are periodic for $n > \max Q$. 
\end{theorem}

\begin{proof}
This is the special case of Theorem~\ref{thm:Qfinite} when $S$ contains just a single positive rod.
\end{proof}

Expansions to a single antirod imply periodicity too.
 
\begin{theorem}\label{thm:periodictoantirod}
    If  $R$ expands to $[\anti{k}]$ via  finite $Q$
    then $R$ expands to $[2k]$ via $Q \cup [\anti{k}] \cup Q\anti{k}$, and
    $\tcount(n,R)$ is eventually periodic with period $2k$.
\end{theorem}
\begin{proof}  
Since $[\anti{k}]$ expands via  $[\anti{k}]$ to $[2k]$, Corollary~\ref{cor:composition} implies $R$ expands to $[2k]$ via $Q \cup [\anti{k}] \cup Q\anti{k}$.
\end{proof}

For example, applying odd sign swap (Corollary~\ref{cor:oddsignswap}) to the expansion in Equation~\ref{eq:period3expansion}
leads to
\begin{equation*} 
    [1, \anti{2}] \trailright{[1]} [\anti{3}].
\end{equation*}
Then the previous theorem implies
\begin{equation}\label{eq:period6}
    [1, \anti{2}] \trailright{[1, \anti{3},\anti{4}]} [6].
\end{equation}
We see one more time that the train counts for
 $[1,\anti{2}]$ start $1,1,0,-1,-1,0$ and then repeat with period $6$, as shown in Table~\ref{table:1anti2} and Equation~\ref{eq:1anti2-6}.  

There are many rod sets $R$ with eventually periodic net train counts, since every sequence of integers counts the net trains for some rod set $R$.
But such $R$ are rarely finite. 

The following characterization of periodic finite rod sets includes a proof that eventually periodic finite rod sets are actually periodic.

\begin{theorem}\label{thm:periodicityequivalents}
The following assertions about finite rod set $R$ are equivalent.

\begin{enumerate}
    \item The net train counts for $R$ are eventually periodic.
    \item     $\tcount(n,R)$ is bounded.
    \item     $\tcount(n,R)$ contains two identical sequences of length $\max R$.
    \item The net train counts $\tcount(n,R)$ are periodic starting at $n=0$.
\item The polynomial $1 -   \rodgeneratingfunction(x,R)$ divides $1- x^p$ for some $p$.
\item $1 - \rodgeneratingfunction(x,R)$ is a product of distinct cyclotomic polynomials.

\item  The recursion generating function   for  $R$ is a finite sum of
  geometric   series of the form
  \begin{align}\label{eq:periodicgeneratingfunction}
  \sum_{n=0}^\infty \tcount(n, R)x^n
    &=  \frac{q(x)}{1-x^m}\\ \notag
    &= q(x)\left(1 + {x^m} +{x^{2m}} + \cdots \right),
  \end{align}
  where
  \begin{equation*}
    q(x) = \sum_{k=0}^{m-1}\tcount(k,R)x^k \ .
\end{equation*}

\end{enumerate}
\end{theorem}
\begin{proof}
Clearly (1) implies (2)
 even without the finiteness hypothesis.

The pigeonhole principle shows that (2) implies (3): there are only finitely many bounded sequences of length $\max R$. 

To show (3) implies (4), we start with two identical  sequences of length $m=\max R$ starting at $a$ and $a+p$.

The recursion shows that $\tcount(a+m)=\tcount(a+p+m)$ and then by induction that  for all $n \ge a$, $\tcount(n)=\tcount(n+p)$. Thus $\tcount(n)$ is eventually periodic with period $p$. 

Now we need to show that $\tcount(n)$ is also periodic for $n<a$, which we can do by ``running the recurrence relation backwards" to compute earlier values of $n$. We start with
\begin{align*}
    \tcount(a+m-1) &= \sum_{k \in R} \tcount(a+m-1-k) \\
    &= \tcount(a-1) + \sum_{k \in R, k<m} \tcount(a+m-1-k),
\end{align*}
and compute
\begin{align*}
    \tcount(a-1) &= \tcount(a+m-1) -  \sum_{k \in R, k<m} \tcount(a+m-1-k) \\
    &=\tcount(a+p+m-1) -  \sum_{k \in R, k<m} \tcount(a+p+m-1-k)\\
    &=\tcount(a+p-1).
\end{align*}

We continue by induction to compute earlier terms and conclude that starting at $n=0$, $\tcount(n,R)$ is periodic with period $p$.
Clearly (4) implies (1), so the first four assertions are equivalent.

The last three equivalences involve the rod count and recursion generating functions for $R$.

Theorem~\ref{thm:periodicity}  says periodicity is equivalent to $R$ expanding to the one element rod set $[p]$, for which $\rodgeneratingfunction(x,[p]) = x^p$.
When $R$ is finite 
Theorem~\ref{thm:productformula} then implies (1) and (5) are equivalent.

The equivalence of (5) and (6) follows from the fact that 
the complete factorization of $x^p-1$ into irreducible polynomials is
\begin{equation}\label{eq:palindromy}
    x^p - 1 = \prod_{d | m} \Phi_d(m),
\end{equation}
where $\Phi_d(x)$ is the cyclotomic polynomial whose
roots are the primitive $d$th roots of unity\cite{cyclotomy}.

The equivalence of (7) is Theorem~\ref{thm:generatingfunction}.
\end{proof}

Repeated cyclotomic factors lead to a rod set that isn't periodic. For example, if
$1 - \rodgeneratingfunction(x,R) = (1-x)^2 = 1-2x+x^2$ then $R = [1,1,\anti{2}]$. The expansion
\begin{equation*}
[1] \trailright{[\anti{1}]}[1,1,\anti{2}] 
\end{equation*}
has dual
\begin{equation*}
    [1,1,\anti{2}]\trailright{\rod{\trains[1]}}[1]. 
\end{equation*}
Here $Q = \rod{\trains[1]}$ is infinite and the
net train count sequence
 $
\tcount(n, [1,1,\anti{2}]) = n+1
$
is unbounded.

In Section \ref{sec:duality} on Duality, we discussed examples where the  rod lengths were infinite arithmetic progressions. Some rod sets where the lengths are finite arithmetic progressions also provide examples of periodicity.

The expansion
\begin{equation*}
   [\anti{1}, \anti{2}, \ldots, \anti{m-1}] 
   \trailright{[\anti{1}]}    [m],
\end{equation*}
corresponding to the factorization
\begin{equation*}
    1-x^m = (1-x)(1 + x + x^2 + \cdots + x^{m-1}),
\end{equation*}
tells us that
\begin{align*}
  \tcount(n, [\anti{1}, \anti{2}, \ldots, \anti{m-1}]) &=
  \begin{cases}
    1 & \text{ if } n \equiv 0 \pmod{m} \\
  -1 & \text{ if } n \equiv 1 \pmod{m} \\
    0 & \text{ otherwise}.
    \end{cases} 
\end{align*}
Similar results follow from the expansions
\begin{align*}
     [1, \anti{2}, 3, \anti{4}, \ldots, \anti{2m}]
  & \trailright{[1]}    [\anti{2m+1}] \\
  \text{and }
    [1, \anti{2}, 3, \anti{4}, \ldots, 2m-1]
   & \trailright{[1]}    [2m].
\end{align*}

When $m$ has multiple prime factors the polynomial $1-x^m$ can lead to interesting periodic rod sets. 
 For example, 
\begin{equation*}
    [\anti{1},3,4,5,\anti{7}, \anti{8}] \trailright{} [30].
\end{equation*}
We encourage the reader to check this result using 
Algorithm~\ref{alg:findQ}.
Theorem~\ref{thm:periodictoantirod} says we can stop when we've expanded  to $\anti{15}$. 
This expansion corresponds to 
rewriting Equation~\ref{eq:palindromy} as
\begin{align*}
    1-x^{30} &= (1+x-x^3-x^4-x^5+x^7+x^8)\\
    &\quad \times (1 - x + x^2 + x^5 - x^6 + x^7 - x^{15} + x^{16} - x^{17} - x^{20} + x^{21} - x^{22}) \\
    &=  \Phi_{30}(x) \times (-1)\Phi_{1}(x)\Phi_{2}(x)\Phi_{3}(x)\Phi_{5}(x)\Phi_{6}(x)\Phi_{10}(x)\Phi_{15}(x).
\end{align*}

In the examples of periodic rod sets that we've seen so far, there has been at most one rod of each length, but there are periodic rod sets $R$ that contain multiple rods of the same length. One such $R$ comes from factoring the cyclotomic polynomial $\Phi(105)$, the first cyclotomic polynomial with a nonzero coefficient different from $\pm 1$. The rod set $R$ has period
 $105$ with $\max R=48=\varphi(105)$ and two rods each of lengths $7$ and $41$. 

 \section{Expandability and Scaling}

In this section we study expansions 
of $R$ to $S$ via $Q$ when all three rod sets are finite and we specify only the lengths of the rods in $S$, not their multiplicities.  
After some general definitions we focus on 
expansions to one or two rod lengths. Our theorems depend on interesting divisibility conditions on the net train counts of $R$. They lead to more identities, to combinatorial proofs about divisibility in Fibonacci and Lucas sequences, and to a combinatorial view of quadratic factors of Borwein trinomial.

\begin{definition}
    The expansion \RQS{} is a \emph{finite} expansion when all three rod sets $R$, $Q$ and $S$ are finite.
\end{definition}
 
As usual, we will write
 $\tcount(n)$ for $\tcount(n,R)$ when there is no potential ambiguity.  

\begin{definition}
A \emph{shape} is a set of distinct rod lengths, usually written as an ordered list in angle brackets. When shape $H$ is finite we write $|H|$ for its cardinality.
     The shape of rod set $R$, $\shape (R)$, is the set   of lengths of rods in the reduced rod set $R'$ equivalent to $R$; if $R$ is finite, so is its shape.  For each length $\ell \in \shape(R)$, the \emph{multiplicity}of $\ell$ is 
     $\rodcount(\ell,R') \ne 0$, the net count of rods of that length in $R'$. The rods in $R'$ of length $\ell$ are antirods just when the multiplicity is negative.
\end{definition}
We can rewrite the recursion for $R$ using shape and multiplicities:
\begin{align}\label{eq:recursionwithshape}
    \tcount(n,R) & = \ \quad \sum_{r \in R}\  (\sign r)F(n-r,R)\notag\\ 
    &= \sum_{\ell \in \shape(R)} \rodcount(\ell, R)\tcount(n-\ell,R).
\end{align}

For example, if $R=[1, \anti{1}, 2,2,2,\anti{3}, \anti{3}, 3] \equiv [2^3, \anti{3}^1]$, then $\shape (R)=\langle 2,3 \rangle$ and the multiplicities of 2 and 3 are 3 and $-1$, respectively.

\begin{definition}
If there is a finite expansion of rod set $R$ to a rod set with shape $H$ then $H$ is necessarily finite with cardinality $n$. Then we say that $R$  is \emph{$n$-expandable} and \emph{expands to} $H$.
\end{definition}

Note that any rod set expands trivially to its shape via the empty $Q$. 

Our goal is to characterize finite expansions of a rod set $R$ to shapes $H$ of cardinality 1 and 2 in terms of conditions on the net train counts
$\tcount(n,R)$.

For example, the expansion  
       \begin{equation*}
          R = [\anti{1},\anti{1}, \anti{2},\anti{2}]   \trailright{[\anti{1},\anti{1},2.2]}
           [\anti{4}^4] = S,
       \end{equation*}
       corresponding to the polynomial factorization
       \begin{equation*}
           1+4x^4 = (1-2x + 2x^2)(1+2x + 2x^2),
       \end{equation*}
shows that $[\anti{1},\anti{1}, \anti{2},\anti{2}] $ is 1-expandable to $\langle 4 \rangle$.
The train counts are
\begin{equation}\label{eq:trainsanti1^2anti2^2}
       \tcount(n,R) =  1, -2, 2, 0, -4, 8, -8, 0, 16, -32, 32, 0, -64, 128, -128, 0 \ldots.
       \end{equation}
Since $\max Q=2$, Theorem~\ref{thm:Qfinite}, tells us that the discrepancy $\discrepancySeq(n,R, [\anti{4}^4])$  must vanish
for $n \ge 3$. In particular, 
$\tcount(3)+4F(3-4) = \tcount(3) = 0$ and
 $\tcount(4)+4F(0)=  -4 + 4 = 0$.

If $\max S - \max R = \max Q$ were greater, then we'd have more discrepancies that would need to be $0$. 
The next theorem generalizes these observations.

\begin{theorem} \label{thm:1elementshape}
There is a finite expansion of rod set $R$ to a rod set $S$ with shape $\langle a \rangle$ if and only if $ \tcount(a-i,R) = 0$ for   $i = 1, 2, \ldots, \max R - 1$ and $\tcount(a,R) \ne 0$. Then $\rodcount(a,S) = \tcount(a,R)$.
\end{theorem}
\begin{proof}
Suppose $\shape(S) = \langle a \rangle $ and $Q$ is finite in the expansion $R \trailright{Q} S$. 

Theorem~\ref{thm:Qfinite} implies $a=\max R + \max Q$, so the discrepancies measuring the ways the counts $\tcount(n,R)$ fail to satisfy the recursion for $S$ are $0$ for $n> \max Q  = a - \max R$. Thus for positive $i <\max R$ 
    \begin{align*} 
    \tcount(a-i,R) &= \tcount(a-i,R)
         - \rodcount(a,S)\tcount(a-i-a,R) \\
        &= \discrepancySeq(a-i,R,S) \\
         &= 0 .
    \end{align*}     
Now consider the discrepancy at $a$: 
     \begin{equation}\label{eq:solvemultiplicity}
        \discrepancySeq(a,R,S)=\tcount(a,R)-\rodcount(a,S)\tcount(a-a,R)=0.
        \end{equation}
Since $\tcount(a-a,R) = 1$  we have $\rodcount(a,S) = \tcount(a,R) \ne 0$. 

To prove the converse, suppose $ \tcount(a-i,R) = 0$ for   $i = 1, 2, \ldots, \max R - 1$.
Let $Q$ be the rod set that expands $R$ to the rod set $S$ with $\tcount(a,R)$ copies of $a$, with $a$ an antirod if $\tcount(a,R)$ is negative. 
 The assumptions on the net train counts of $R$ say that
 the discrepancy $\rodcount(n,Q) = \discrepancySeq(n,R,S) = 0$  when $a - \max R < n < a$. Setting $\rodcount(a,S) = \tcount(a,R)$ guarantees discrepancy $0$ at $n=a$.
 Then 
 Corollary~\ref{cor:maxR0inarow} implies $Q$ is finite.
\end{proof}

Next we apply similar discrepancy techniques to the more complex characterization of $2$-expansions. Finite expansions of a finite $R$ to $S$ with shape $\langle a,b \rangle$  correspond to recursion identities of the form
\begin{equation*}
    \tcount(n,R)=  \rodcount(a,S)\tcount(n-a,R) + \rodcount(b,S)\tcount(n-b,R) 
\end{equation*}
and to trinomial multiples of 
the polynomial $1-\rodgeneratingfunction(x,R)$.

 For example there is a finite expansion of the Fibonacci rod set $[1,2]$ to $[8^{48},\anti{12}^{7}]$,
corresponding to the Fibonacci identity
\begin{equation*}
    \tcount(n) = 48\tcount(n-8) - 7 \tcount(n-12),
\end{equation*}
for large enough $n$, and to the factorization
\begin{equation*}
    1-48x^{8} + 7x^{12}= q(x)(1-x-x^2),
\end{equation*}
where
\begin{equation*}
  q(x)=  1 + x + 2 x^2 + 3 x^3 + 5 x^4 + 8 x^5 + 13 x^6 + 21 x^7 - 14 x^8 + 7 x^9 - 7 x^{10}.
\end{equation*}

In the next section we will show how our methods lead to similar Fibonacci identities and discuss how they appear in the literature.  

In 
Theorem~\ref{thm:1elementshape} we saw that 1-expansions of $R$ correspond to sequences of $0$'s in $\tcount(n,R)$.  For 2-expandability the correspondence is to divisibility among net train counts. For example, we will see soon that the expansion above of $[1,2]$ to $[8^{48},\anti{12}^{7}]$ is equivalent to the fact that
\begin{equation*}
  48\tcount(3,[1,2]) = 48 \times 3 = 144 = \tcount(11,[1,2]).
\end{equation*}
The Padovan expansion of $[2,3]$ to $[4^3,13]$ in 
Figure~\ref{fig:handcalculation} corresponds to the adjacent pair of scalings by a factor of 3:
\begin{align*}
      3\tcount(7,[2,3]) &= 3 \times 3 = \ 9= \tcount(11,[2,3])\\
         3\tcount(8,[2,3]) &= 3 \times 4 = 12= \tcount(12,[2,3])
\end{align*}
in the Padovan sequence
\begin{equation*}
    P(n) =1, 0, 1, 1, 1, 2, 2, \textbf{3}, \textbf{4}, 5, 7, \textbf{9}, \textbf{12},16, 21,28,  \ldots 
\end{equation*}

To prove our general theorem on $2$-expandability, we formally define scaling.

\begin{definition}\label{def:scaling}
    The finite rod set $R$ \emph{scales at counts $(c,b)$} when there is a unique integer \aaa{}  such that
    \begin{equation}\label{eq:scaletraincounts}
       \tcount(b-i, R) = \aaa \tcount(c-i, R) 
    \end{equation}
    for $1 \le i < \max R$.
    Then we say $R$  scales by $\aaa$ at $(c,b)$.
 
\end{definition}

The requirement that $\alpha$ be unique excludes the degenerate case where all the Equations~\ref{eq:scaletraincounts} are of the form $0=\alpha \times 0$. 

The previous examples show that $[1,2]$ scales by $48$ at $(4,12)$ and 
$[2,3]$ scales by $3$ at $(9,13)$. 

For a larger example including antirods we revisit the Padovan expansion 
\begin{equation*}
    [2,3] 
\trailright{Q = [2,3, \anti{4}^2,5^2,\anti{6},8, \anti{9},10]}
[4^3,13].
\end{equation*}
The Exchange Corollary~\ref{cor:switchRQ}
implies
\begin{equation}
    \anti{Q} \trailright{[\anti{2}, \anti{3}]} [4^3,13].
\end{equation}
The net train counts for $\anti{Q}$ are
\begin{equation*}
    \tcount(n,\anti{Q}) = 1, 0, -1, -1, 3, 0, -3, -3, 9, 0, -9, -9, 27, 1, -27, -28, 80, 6, -81, \ldots .
\end{equation*}
In the following table the 
the third line shifts the net train counts in the second by $4$.
Then the second line  is 3 times the third between $n=4$ and $n=12$.
\begin{equation*}
    \begin{array}{rrrrrrrrrrrrrrrrrr}
 n:&    0& 1& 2& 3& 4& 5& 6& 7& 8& 9& 10& 11& 12& 13&   \ldots\\
\tcount(n):& 1& 0&-1&-1& 3& 0&-3&-3& 9& 0& -9& -9& 27&  1&  \ldots \\
\tcount(n-4):&&  &  &  & 1& 0&-1&-1& 3& 0& -3& -3&  9&  0&   \ldots \\
    \end{array}
\end{equation*}
Therefore $\anti{Q}$ scales by $3$ at $(4,13)$. 

The next theorem shows that scaling is the heart of  $2$-expandability.

\begin{theorem}\label{thm:2elementshape}
There is a finite expansion of rod set $R$ to a rod set $S$ with shape $\langle a,b \rangle$
if and only if there is a nonzero integer \aaa{} such that
       $R$ scales by \aaa{} at $(b-a,b)$, and $\tcount(b,R) \ne \alpha \tcount(b-a,R)$. Then
       \begin{align}
       \rodcount(a,S) &= \aaa \label{eq:na}  , \\ 
  \rodcount(b,S) & =\tcount(b,R)-\aaa\tcount(b-a,R).  \label{eq:nb} 
       \end{align}
\end{theorem}
\begin{proof}
Suppose there is a finite expansion   $\RQS$  with $\shape(S) = \langle a,b \rangle$. 
Theorem~\ref{thm:Qfinite} implies $b=\max R + \max Q$, so the discrepancies measuring the ways the counts $\tcount(n,R)$ fail to satisfy the recursion for $S$ are $0$ for $n>b - \max R$. In the range $b - \max R < n < b$, we have $\tcount(n-b,R)=0$.   

At $n=b-i$ for positive $i <\max R$ 
we compute
   \begin{align*}
         \discrepancySeq(b-i,R,S)&=\tcount(b-i,R)-\rodcount(a,S)\tcount(b-a-i) 
         - \rodcount(b,S)\tcount(b-b-i) \\
         &= \tcount(b-i,R)-\rodcount(a,S)\tcount(b-a-i) \\
         &= 0 .
    \end{align*}
Thus for   $i = 1, 2, \ldots, \max R - 1$
\begin{equation} 
     \tcount(b-i,R) = \rodcount(a,S)\tcount(b-a-i,R),
\end{equation}
so $R$ scales by $\aaa  = \rodcount(a,S) $ at $(b-a,b)$.

Now consider the discrepancy at $b$: 
     \begin{equation}\label{eq:solvenb}
        \discrepancySeq(b,R,S)=\tcount(b,R)- \rodcount(a,S)\tcount(b-a)-\rodcount(b,S)\tcount(b-b,R)=0.
        \end{equation}
Since $\tcount(b-b,R) = 1$  we can solve for $\rodcount(b,S)$ to reach Equation~\ref{eq:nb}. Since $\tcount(b,R) \ne \alpha \tcount(b-a,R)$, $\rodcount(b,S) \ne 0$. When $\rodcount(b,S)$ is negative, $b$ is an antirod.

To prove the converse, suppose that  
$R$ scales by some \aaa{} at $(b-a,b)$.
Use Equations~\ref{eq:na} and \ref{eq:nb} to define $\rodcount(a,S)$ and $\rodcount(b,S)$. 

Let $Q$ be the rod set that expands $R$ to the $S$ with shape  $\langle a,b \rangle$ and those net rod counts.
 It follows that the discrepancy $\rodcount(n,Q) = \discrepancySeq(n,R,S) = 0$  when $b - \max R < n \le b$. 
 Corollary~\ref{cor:maxR0inarow} implies $Q$ is finite.
\end{proof}
 When a rod set is $1$- or $2$-expandable there is at most one $S$ for each possible such expansion.
 The existence theorems above determine the multiplicities in $S$. For $3$-expandability and higher it's easy to find counterexamples to uniqueness. For example, 
 
\begin{equation*}
    [1,2,2]\trailright{[1,2]}[2^2,3^3,4^2]
\end{equation*}
and 
 \begin{equation*}
    [1,2,2]\trailright{[1,2,2]}[2,3^4,4^4]
\end{equation*}
have the same shape, $\langle 2,3,4 \rangle$.

\section{Fibonacci and Lucas 2-Expansions}

In this section we study $2$-expansions of rod sets with shape $\langle 1,2 \rangle$ to prove identities and divisibility properties for their train counts.
We first introduce new notation that will greatly simplify our work. 

\begin{definition}
Suppose rod set $R$ contains $k$ copies of rod $r$. Then we list its contents in double square brackets as
\begin{equation*}
    R = [[ \ldots, \length(r)^{\sign(r)k}, \ldots]].
\end{equation*}
For example
\begin{equation*}
[1^2,\anti{2}^3] = [[1^2,2^{-3}]].
\end{equation*}
\end{definition}

Typographically, this just slides the overbar signaling an antirod right to become a minus sign in the exponent counting the multiplicity. 

Now collecting copies of rods of the same length becomes adding exponents, with the \emph{caveat} that $r^0$ means no rods of length $r$, not $1$. For example
\begin{equation*}
    [1^2, 2^5, \anti{2}^8,3, \anti{3}]
    =    [[1^2, 2^5, 2^{-8},3, 3^{-1}]]
    \equiv [[1^2,2^{-3}]].
\end{equation*}

Expansions work nicely too. Expanding both $1$'s in $[[1^2,2^{-3}]] $ yields
\begin{equation*}
     [[1^2,2^{-3}]] \trailright{ [[1^2]]} 
     [[(1+1)^{2 \times 2},(2+1)^{-3 \times 2}, 2^{-3}]] 
     \equiv [[2,3^{-6}]].
\end{equation*}

\begin{definition}\label{def:lucasrodset}
A \emph{Lucas} rod set $R$ is a rod set of the form $[[1^s,2^t]]$,
where $s$ and $t$ are nonzero, 
$\gcd(s,t) = 1$,  and
$\tcount(n,R) \ne 0$ for $n \ge 0$. 
\end{definition}

The double bracket notation for rod sets makes writing recursions easy. The exponents are the coefficients.
The train count sequence for the Lucas rod set $[[1^s,2^t]]$ satisfies the recursion
\begin{equation*}
    \tcount(n, [[1^s,2^t]]) = s \tcount(n-1,[[1^s,2^t]]) + t \tcount(n-2,[[1^s,2^t]])
\end{equation*}
whatever the signs of the rods.
Sequences satisfying that recursion (with the indices shifted by $1$) are known in the literature as regular, nondegenerate fundamental Lucas sequences  \cite{lucas-book}. 

The requirement that no Lucas rod set train count is $0$ is less restrictive than it appears. It is clearly true when $s$ and $t$ are positive. Odd sign swap shows it's true when just $t$ is positive.
When $t$ is negative it fails for the two periodic rod sets
$[[1,2^{-1}]]$ and $[[1^{-1},2^{-1}]]$.
The proof of Theorem 38 in \cite{lucas-book} shows using quadratic fields that these are the only two rod sets with shape $\langle 1,2 \rangle$ with a $0$ train count. We have not tried to prove this using our methods.

In Theorem~\ref{thm:lucastotwo} we will systematically construct three classes of  expansions of Lucas rod sets to rod sets with a two element shape. We start with two easy lemmas.

\begin{lemma}\label{lemma:sgreaterthan1}
Let $R$ be a Lucas rod set with $s>1$. Then for $n \ge 1$,  $s \mid \tcount(n,R)$ if and only if $n$ is odd. 
\end{lemma}
\begin{proof}
The $s$ rods of length 1 come in  $s$ colors, so $s$ divides the number of trains of length $n$ that contain a rod of length $1$. Every train of odd length must include at least one rod of length $1$. 

When $n \ge 2$ is even, the trains with no rod of length $1$ are those made from only rods of length $2$; there are $t^{n/2} \ne 0$ of them. Since $s$ divides the number of rods that contain a rod of length $1$, whether or not it divides $\tcount(n)$ depends on whether or not it divides $t^{n/2}$. Since by definition, $\gcd(s,t)=1$, $s$ does not divide $t^{n/2}$.

\end{proof}

\begin{lemma}\label{lemma:Lucasto1}
No Lucas rod set has a finite expansion to
a shape with cardinality $1$.
\end{lemma}
\begin{proof}
If rod set $R$ has a finite expansion to $\langle c \rangle$ then  Theorem~\ref{thm:1elementshape} implies \\ $\tcount(c-1, R)$ = 0. That never happens for a Lucas rod set.
\end{proof}

We will illustrate the proof of the next theorem using the running example $[[1^2,2^{-3}]]$ with
net train counts starting at $\tcount(0) =1$
\begin{equation}\label{eq:Lucas1^3anti2^2}   1,2,1,-4,-11,-10,13,56,73,-22,-263,-460,-131,1118,2629,1904, \ldots\ \footnote{For more on this sequence, see OEIS A088137}  .
\end{equation}

Lemma~\ref{lemma:sgreaterthan1} explains why every other count is even. Notice too that so far every third count is a multiple of $1$, every fourth a multiple of $4$, and every fifth a multiple of $11$. At the end of this section we will prove that this pattern continues. 

\begin{theorem}\label{thm:lucastotwo}
For each of the following pairs $a < b$ of positive integers, the Lucas rod set $R$ expands via a finite rod set to a unique rod set with shape $\langle a,b \rangle$.  
\begin{enumerate}

\item $(a,a+1)$.

\item $(a, a+2)$, when $s=1$ or when $a$ is even.

\item $(kd, (k+1)d)$ for $d > 2$ and $k \ge 1$.
 
\end{enumerate}
\end{theorem}

\begin{proof}
We need prove only existence; Theorem~\ref{thm:2elementshape} guarantees uniqueness. Since by definition Lucas rod sets have nonzero train counts, Lemma~\ref{lemma:Lucasto1} guarantees that the expansions below never degenerate to shapes of cardinality $1$.

\begin{enumerate}
    \item 
       To find the expansion to shape $\langle a,a+1 \rangle$
    we start from $[[1^s,2^t]] = [[1^{\tcount(1)},2^{t\tcount(0)}]]$ and repeatedly expand all of the smallest rods.  The algebra that computes the exponent of the smaller rod is  the same algebra that computes the next term in the recursion. At each step the expanded rod set is
\begin{equation}\label{eq:lucaseq1}
    [[a^{\tcount(a)},(a+1)^{t\tcount(a-1)}]].
\end{equation}

    In the example, we have 
    \begin{equation*}
        [[1^2,2^{-3}]]\trailright{}[[2,3^{-6}]]\trailright{}[[3^{-4},4^{-3}]]\trailright{}[[4^{-11},5^{12}]]\trailright{} 
[[5^{-10},6^{33}]]\trailright{}\ldots.
    \end{equation*}
    
 \item
    To find the expansion to shape $\langle a,a+2 \rangle$ when $s=1$ or $a$ is even, we note that then Lemma~\ref{lemma:sgreaterthan1} implies  $t\tcount(a-1)/s$ is an integer.
    We start with any rod set of the form in Equation~\ref{eq:lucaseq1} and
add $t\tcount(a-1)/s$ rod/antirod pairs of length $a$, expand the $a^{-t\tcount(a-1)/s}$ rods,  and simplify. 
\begin{align}\label{eq:lucaseq2}
          [[a^{F(a)},(a+1)^{tF(a-1)}]] &\equiv [[a^{F(a)+t\tcount(a-1)/s},a^{-t\tcount(a-1)/s}, (a+1)^{tF(a-1)}]] \notag \\&\trailright{[[a^{-t\tcount(a-1)/s}]]} 
          [[a^{F(a)+t\tcount(a-1)/s},(a+2)^{-t^2F(a-1)/s}]] \notag \\
          &= [[a^{\tcount(a+1)/s},(a+2)^{-t^2F(a-1)/s}]] .
      \end{align}

In the example we have
      \begin{equation*}
          [[4^{-11},5^{12}]] \equiv [[4^{-11},4^{6},4^{-6}, 5^{12}]] \trailright{[[4^{-6}]]}[[4^{-5},5^{-12},5^{12},6^{18}]]\equiv [[4^{-5},6^{18}]].
      \end{equation*}

\item To find the expansion to shape $\langle kd, (k+1)d \rangle$ for $d > 2$ and $k \ge 1$ we first explain how to find examples of shape $\langle d,2d \rangle$.

We start with a rod set 
from Equation~\ref{eq:lucaseq1} and add some rod/antirod pairs:
 \begin{align*}
   [[d^{\tcount(d)},(d+1)^{t\tcount(d-1)}]] 
    &\equiv [[d^{\tcount(d)},d^{t\tcount(d-2)},d^{-t\tcount(d-2)},(d+1)^{t\tcount(d-1)}]] \\
    &= [[d^{\tcount(d)+t\tcount(d-2)},d^{-t\tcount(d-2)},(d+1)^{t\tcount(d-1)}]] \\
        &= [[d^{c},d^{-t\tcount(d-2)},(d+1)^{t\tcount(d-1)}]] ,
 \end{align*}
 writing $c=\tcount(d)+t\tcount(d-2)$. Then we expand the rods $d^{-t\tcount(d-2)}$:
\begin{align*}
    [[d^c,d^{-t\tcount(d-2)},(d+1)^{t\tcount(d-1)}]] &\trailright{}  
    [[d^c,(d+1)^{t(\tcount(d-1)-s\tcount(d-2))},(d+2)^{-t^2\tcount(d-2)}]] \\ 
    &\phantom{a} \equiv [[d^c,(d+1)^{t^2\tcount(d-3)},(d+2)^{-t^2\tcount(d-2)}]].
\end{align*}
We continue expanding the rods of intermediate length $d+i$ to produce rod sets of shape $\langle d, d+i, d+i+1 \rangle$  until $i=d-2$:
\begin{align}\label{eq:lucaseq3}
    [[d^{\tcount(d)},(d+1)^{t\tcount(d-1)}]]  &\trailright{} [[d^c,(d+1)^{t^2\tcount(d-3)},(d+2)^{-t^2\tcount(d-2)}]]   \notag\\ 
    &\trailright{} [[d^c,(d+2)^{-t^3\tcount(d-4)},(d+3)^{t^3\tcount(d-3)}]] \notag \\
    &\trailright{}[[d^c,(d+3)^{t^4\tcount(d-5)},(d+4)^{-t^4\tcount(d-4)}]] \notag\\
    &\ldots \notag\\
    &\trailright{} [d^{c},(2d)^{(-1)^{d-1}t^d}].
\end{align}

Here is the expansion for $d=4$ in our example,
\begin{align*}
   [[4^{-11},5^{12}]] &\equiv [[4^{-11},4^{-3},4^{3},5^{12}]]\\
   &\trailright{[[4^{3}]]}[[4^{-14},5^{18},6^{-9}]] \\
   &\trailright{[[5^{18}]]}[[4^{-14},6^{27},7^{-54}]]\\
   &\trailright{[[6^{27}]]}[[4^{-14},8^{-81}]].
\end{align*}

To find the rest of the expansions of $R$ to shapes $\langle kd, (k+1)d \rangle$, we  start with the rod set $[[d^{c},(2d)^{(-1)^{d-1}t^d}]]$ with shape $\langle d, 2d \rangle$ and repeatedly expand all of the smallest rods to obtain a rod set of shape $\langle 2d, 3d \rangle$. We proceed inductively, expanding all of the rods of length $kd$ in the rod set of shape $\langle kd, (k+1)d \rangle$ to obtain a shape that increases $k$ by $1$. Corollary~\ref{cor:composition} implies that the expansions we create are also finite expansions of $R$. 

Below are first steps in the expansion,
\begin{align}\label{eq:lucaseq4}
    [[d^{c},(2d)^{(-1)^{d-1}t^d}]] & \trailright{[[d^{c}]]} [[(2d)^{c^2+(-1)^{d-1}t^d},(3d)^{(-1)^{d-1}ct^d}]] \notag \\
     & \trailright{[[(2d)^{c^2+(-1)^{d-1}t^d}]]} \ldots
\end{align}

The formulas for the multiplicities get increasingly complicated, but fortunately, we don't need their general expressions. Here are expansions starting from shape  $\langle 4, 8 \rangle$ to shapes $\langle 8,12 \rangle$ and $\langle 12,16 \rangle$.
\begin{align*}
    [[4^{-14},8^{-81}]] & \trailright{[[4^{-14}]]} [[8^{(-14)(-14)-81},12^{(-14)(-81)}]] \\
    & = \phantom{aaaaaa} [[8^{115},12^{1134}]]\\
    & \trailright{[[8^{115}]]}
    [[12^{(115)(-14)+1134},16^{(115)(-81)}]] \\
    & = \phantom{aaaaaa} [[12^{-476},16^{-9315}]].
\end{align*}
 \end{enumerate}
\end{proof}

The expansions we found in Theorem~\ref{thm:lucastotwo} correspond to known identities for Lucas rod set train counts. For example,
Equation~\ref{eq:lucaseq1} corresponds to the identity 
\begin{equation*}
       F(n)= \tcount(a)F(n-a)+t\tcount(a-1)F(n-a-1),
    \end{equation*}
which is quadratic, because the multiplicities in the expansion are themselves rod counts.
This identity is a more general version of Identity $3$ in 
Benjamin and Quinn~\cite{benjamin}. The combinatorial proof there generalizes to our Lucas rod sequence identity, with the same indexing.

Reindexed, the recursions associated with Equations~\ref{eq:lucaseq2}, \ref{eq:lucaseq3}, and \ref{eq:lucaseq4} can be found as identities derived from characteristic roots. The recursion from Equation~\ref{eq:lucaseq2} follows in a few steps from Equation $2.43$ in Ballot and Williams~\cite{lucas-book}. The recursion from Equation~\ref{eq:lucaseq3} is Equation $2.16$ there. When $R$ is the Fibonaccis, the exponent $c$ in  Equation~\ref{eq:lucaseq3} is a Lucas number, and the corresponding recursion is Identities $15$a and $15$b in Vajda~\cite{vajda}. The recursions from Equation~\ref{eq:lucaseq4} are related to Equations $2.16$ and $2.20$ in Ballot and Williams~\cite{lucas-book}. To reach our results from theirs we need to use divisibility properties that we establish in the next theorem.

\begin{definition}
    Let $R$ be a rod set. The sequence $\trains(n,R)$ is a \emph{divisibility sequence} when
    \begin{equation*}
    m+1 \mid n+1 \implies \tcount(m) \mid \tcount(n).
        \end{equation*}
\end{definition}

Note that this definition would be nicer and match what most poeple use in the literature if we removed $\tcount(0)= 1$ from the sequence and reindexed so that the displayed implication was  $m \mid n \implies \tcount(m) \mid \tcount(n)$.

\begin{theorem}\label{thm:lucasdivseq}
The net train counts $\tcount(n,R)$ for a Lucas rod set $R$ form a divisibility sequence.
\end{theorem}

\begin{proof}
Suppose $m+1 \mid n+1$. 

If $m=0$ then $m+1 =1= \tcount(0)$ divides  $\tcount(n)$ for all $n$. 

If $m=1$ then $\tcount(m)=\tcount(1)=s$. If $s=1$  then $s \mid \tcount(n)$ for all $n$. If $s > 1$ then $m+1 = 2 $ divides $n+1$ so $n$ is odd, and Lemma~\ref{lemma:sgreaterthan1}
 shows that
$\tcount(1) = s$ divides $\tcount(n)$.

If $m > 1$ we use
Theorem~\ref{thm:2elementshape}, which shows that if rod set $R$ has a finite expansion to shape 
$\langle n-m, n+1 \rangle$ then there is a nonzero $\alpha$ such that $\tcount(n) = \alpha\tcount(m)$, so $\tcount(m) | \tcount(n)$.
We show there is such an expansion.
We know that for some $k$, 
$n+1 = k(m+1)$.
Then
\begin{align*}
    \langle n-m,\  n+1 \rangle
        &=   \langle k(m+1)-1-m,\  k(m+1) \rangle \\
        &=   \langle (k-1)(m+1),\  k(m+1) \rangle.
\end{align*}
The lengths in that shape are consecutive multiples of $m+1$ so
 Case (3) of Theorem~\ref{thm:lucastotwo} shows that the desired expansion exists.
\end{proof}

Note that the preceding proof used Theorem~\ref{thm:2elementshape} and then the third case in Theorem~\ref{thm:lucastotwo} only when $m >1$. We could have dealt with $m=0$ and $m=1$ using the first two cases to prove the existence of expansions that imply the divisibility, but the arguments we chose instead are simpler.

The regular, nondegenerate fundamental Lucas sequence $U(P,Q)$ corresponds to the Lucas rod set $R=[[1^P,2^{-Q}]]$ with $U_n(P,Q)=\tcount(n-1,R)$. Leaving the rod set context, we have  shown that regular, nondegenerate fundamental Lucas sequences are divisibility sequences:

\begin{corollary}
    For a regular, nondegenerate fundamental Lucas sequence $U(P,Q)$, if $m \mid n$ then $U_m(P,Q) \mid U_n(P,Q)$.
\end{corollary}

\section{Borwein trinomials}\label{sec:Borwein}

In this section we apply
Theorem~\ref{thm:2elementshape}   to study expansions of the periodic rod sets $[1,\anti{2}]$ and $[\anti{1},\anti{2}]$
--- the rod sets of shape $\langle 1, 2 \rangle$ that are \emph{not} Lucas rod sets. We find all their expansions to
 a two element shape.  These expansions provide another view of results in the literature on factoring  Borwein trinomials.

\begin{definition}
    A \emph{Borwein polynomial} is a polynomial all of whose nonzero coefficients are $\pm 1$. We will assume in addition that the constant term is 1.
\end{definition}

Borwein polynomials correspond to finite rod sets with just one rod of each length.

It is known that when $a$ and $b$ are relatively prime the Borwein trinomial $1 \pm x^a \pm x^b$ is either irreducible or the product of an irreducible and $1 -x +x^2$ or $1 +x +x^2$ (\cite{koley}, \cite{koley2}). The proofs are analytic. In this section we find all  Borwein trinomial multiples of $1  \pm x +x^2$ by finding all possible expansions of the periodic rod sets $[1,\anti{2}]$ and $[\anti{1},\anti{2}]$ to two-rod sets. We do not prove the irreducibility assertions.

\begin{figure}[h]
    \centering
    \includegraphics[width=0.4\textwidth]{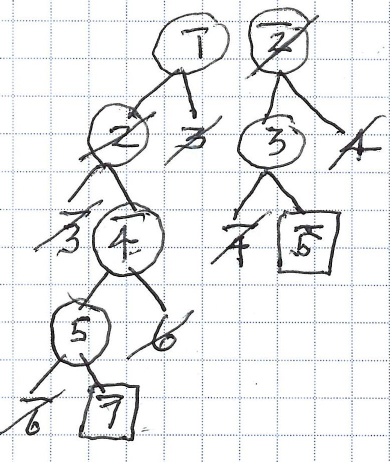}
    \caption{A Borwein expansion}
    \label{fig:borwein}
\end{figure}

Figure~\ref{fig:borwein} shows an expansion tree for the expansion
\begin{equation*}
    [\anti{1}, \anti{2}]
    \trailright{[\anti{1},3,\anti{4},5]}
    [\anti{5},\anti{7}]
\end{equation*}
that corresponds to the factorization
\begin{equation*}
1 + x^5 + x^7 = 
      \left(1 + x  + x^{2}\right) 
      \left(1 - x + x^{3} - x^{4}  + x^{5} \right) .
\end{equation*}

Note that the signs toggle as we move down a level in the tree, and all nodes at the same level have the same sign. 

\begin{theorem}\label{thm:borwein}
   The rod set $[\anti{1},\anti{2}]$ expands by a finite $Q$ to a rod set $S$ with shape $\langle a,b \rangle$ if and only if $a \equiv 1 \pmod{3}$ and $b \equiv 2 \pmod{3}$ (either may be larger). When such an expansion exists, the rods in $S$ are both antirods with multiplicity $1$.

      The rod set $[1,\anti{2}]$ expands by a finite $Q$ to a rod set $S$ with shape $\langle a,b \rangle$ under the same congruence conditions modulo $3$. The rods in $S$ have multiplicity 1; those of even length are antirods.

      In either case the rods in $Q$ have multiplicity $1$.
\end{theorem}
\begin{proof}
To construct the expansions of $[\anti{1},\anti{2}]$ suppose $j$ and $k$ are positive integers and 
\begin{equation*}
   \langle a, b \rangle = \langle 1+3j,2+3k \rangle.
\end{equation*}
Figure~\ref{fig:borwein} illustrates $j=2, k=1$.
Let

\begin{align}
    Q_L &= [\anti{1}, 2, \anti{4}, 5, 
    \ldots, \anti{3j-2}, 3j-1] \notag \\
    Q_R &= [\anti{2}, 3, \anti{5}, 6, 
    \ldots, \anti{3k-1}, 3k].
\end{align}
 Then
 \begin{align*}
     [\anti{1}, \anti{2}]
     &\trailright{Q_L} 
        [\anti{a}, \anti{2}].\\
     [\anti{1}, \anti{2}]
     &\trailright{Q_R} 
        [\anti{1}, \anti{b}]        
 \end{align*}
 and
 \begin{equation}\label{eq:anti1anti2expand}
     [\anti{1}, \anti{2}]
     \trailright{Q_L \cup Q_R}
        [\anti{a}, \anti{b}]   .
 \end{equation}
Rods of the same length appear in the union
$Q = Q_L \cup Q_R$ only as rod/antirod pairs, so all multiplicities in $Q$ are 1.

Odd sign swap (Corollary~\ref{cor:oddsignswap}) then constructs  expansions of $[1, \anti{2}]$ to shape $\langle a,b \rangle$.
    
    That there are no other expansions to two element shapes follows from Theorem~\ref{thm:2elementshape}. Recall that
    \begin{equation*}
        \tcount(n, [\anti{1},\anti{2}]) = 1, -1, 0, 1, -1, 0. \ldots.
    \end{equation*}
    The only possible scalings match $1$ to $-1$ with  $\alpha = -1$ or $1$ to $1$ with  $\alpha = 1$.  Equation~\ref{eq:anti1anti2expand} finds expansions for all cases with $\alpha=-1$. 
    There are no expansions with $\alpha =1$, because  Theorem~\ref{thm:2elementshape} asserts that the scaling must fail somewhere. If $\alpha=1$, the scaling comes from shifting the sequence of rod counts by a multiple of $3$, which aligns the entire sequence. 
\end{proof}

 \begin{corollary}
    The polynomial $1 + x + x^2$ divides the trinomial
    $p(x) =  1 + \alpha x^a + \beta x^b$ if and only  if 
      $a \equiv 1 \pmod{3}$, $b \equiv 2 \pmod{3}$ and $\alpha = \beta = 1$.

    The polynomial $1 - x + x^2$ divides the trinomial
    $p(x) = 1 + \alpha x^a + \beta x^b$ if and only  if 
        $a \equiv 1 \pmod{3}$, $b \equiv 2 \pmod{3}$, $|\alpha| = |\beta| = 1$ and
the signs of terms with even (odd) exponents are positive (negative).

    In either case both $p(x)$  and the quotient are Borwein polynomials.
\end{corollary}

Note that Equation~\ref{eq:anti1anti2expand} tells us the presence and signs of terms in the quotient polynomial for
$Q = Q_L \cup Q_R$.

\section{$2$-Expandability when $\max R > 2$}\label{sec:MaxRGreater2}

So far we have applied Theorem~\ref{thm:2elementshape} when $R$ contains only rods of lengths 1 and 2.
We end our paper  with examples and conjectures on $2$-expansions of rod sets $R$ with $\max R > 2$ and with some ideas for further work.

 When $R=[1,2,...,k]$ the
the net train counts for $R$ start with consecutive powers of $2$ and scale for a while at the start. A finite variant of Equation~\ref{eq:expand123dots}  leads to the expansion
\begin{equation*}
    [1,2,...,k] \trailright{[\anti{1}]}[1^2,\anti{k+1}].
\end{equation*}
Thus for any finite $k$, there exists at least one $2$-expandable rod set $R$ with $\max R=k$. 

Nonetheless, when 
$\max R \ge 3$ the scalability conditions are  quite restrictive.
There are just three rod sets with two different positive rods of length at most $3$. We've studied the Fibonacci rod set $[1,2]$.
To find 2-expansions of the Padovan rod set, $[2,3]$,
we examine the net train count sequence 
\begin{equation*}
    P(n) =1, 0, 1, 1, 1, 2, 2, 3, 4, 5, 7, 9, 12,16, 21,28,37, 49, 86, 114, 151, 200,  \ldots 
\end{equation*}
looking for scalings: two pairs of consecutive terms $(P(b-a-2),P(b-a-1))$ and $(P(b-2),P(b-1))$, where the second pair is an integer multiple of the first. This happens with the net train counts $1,1,1$ when two  $(1,1)$'=s overlap, which corresponds to the expansion to $[1,5]$. 

Then twice we see pairs of consecutive net train counts
 $(1,1)$ and $(2,2)$.  
Further along we find $(3,4)$ and its multiples $(9,12), (12,16)$, and $(21,28)$. These lead to the following Padovan expansions (and also to their odd sign swaps): 
\begin{align*}
 [2,3] &\trailright{[\anti{1},2]} [1,5] \\
        [2,3] &\trailright{[\anti{2},3,\anti{4}]} [2,2,\anti{7}]=[2^2,\anti{7}] \\
    [2,3] &\trailright{[2,\anti{3},4]} [3^2,7] \\
        [2,3] &\trailright{ } [4^3,13] \\
      [2,3] &\trailright{ } [5^4,14]   \\
            [2,3] &\trailright{ } [7^7,16^2] .
\end{align*}
 Those are all the scalings we could see
simply staring at the sequence. Note that past $\tcount(15)=28$ the numbers seem to be growing in a way that makes further such pairs unlikely, which we confirmed with a computer search up to $n=1000$. 

Similar investigations by hand and by computer up to $n=1000$ yield three finite expansions of the Narayana rod set $[1,3]$, to $[2^2,7]$, $[3^3,8]$, and $[11^{67},33]$. 

When we allow antirods but still require $\max R = 3$ and one kind of each rod we find
 the following expansions (and their odd sign swaps):

 There are $2$-expansions from the rod set $[1,\anti{3}]$ to eight more rod sets: $[\anti{4}, \anti{5}]$, $[\anti{4}^2, \anti{7}]$, $[\anti{5}^2, 7]$, $[9^3, \anti{13}]$, $[9^4, 14]$, $[\anti{13}^4, \anti{14}^3]$, $[\anti{13}^7, \anti{16}^3]$, $[\anti{14}^7, 16^4]$. 

 There are $2$-expansions from the rod set $[\anti{2},3]$ to four more rod sets:  $[\anti{5}^2, 7]$, $[\anti{5}^3, 8]$, $[7^3, \anti{8}^2]$, $[\anti{22}^{67}, 33]$.

We've found only two $2$-expansions with $\max R > 3$:
\begin{equation}\label{eq:reflect1}
    [1,\anti{4}]\trailright{[1,2,3,\anti{5},6,9]} [\anti{6}^3, \anti{13}]
\end{equation}
and
\begin{equation}\label{eq:reflect2}
    [3,\anti{4}]\trailright{[3,\anti{4},6,7,8,9]} [\anti{7}^3, \anti{13}] ,
\end{equation}

\section*{What next?}

\begin{itemize}
     \item  There is more to discover on the connections between expansions of $R$ and the structure of the net train count sequence $\tcount(n,R)$.
     Theorems~\ref{thm:1elementshape} and  \ref{thm:2elementshape} settle those connections for expansions to shapes of cardinality 1 and 2.      
    What are the corresponding theorems for larger shapes? 

\item Can our methods prove more of the known results on Borwein polynomials and Lucas sequences? Can they find new results?

\item Can we prove with our methods that Theorem~\ref{thm:lucastotwo} finds all possible $2$-expansions of a Lucas rod set? Can we prove the the train counts of Lucas rod sets are strong divisibility sequences? Are these two tasks related? 

\item We can run the recursion for a finite rod set $R$ backwards. When $R$ contains just one rod of length $\max R$ the values for negative indices are integers. With some reindexing, we can view them as the train counts for a new rod set, the \emph{reflection} of $R$. For example, $[\anti{1},2]$ is the reflection of $[1,2]$. In Equations~\ref{eq:reflect1} and \ref{eq:reflect2} the corresponding rod sets are reflections of each other, and both expansions are valid. What more can we say about reflection?

\item Much of our rod set algebra works when rod lengths are elements of an ordered group. We use the fact that lengths are positive integers only as exponents in the rod generating functions and when discussing the recursions determined by rod sets.   Is any of this generalization useful or surprising?

\item Our theorems on finite expansions translate to theorems about polynomial multiplication. What we don't have are theorems that say there are no finite expansions to particular rod sets. Such theorems would be rod set algebra proofs of irreducibility. Are there any?

\item
 There is substantial interest in factoring polynomials over finite fields. Can our rod and tree model for polynomial arithmetic be of use?

\item We have completely avoided analysis. In particular, we discovered some early examples of expansion  by looking at the growth rates of train counts. 
Can our rod set model say anything useful  about train count growth rates?

\end{itemize}

\section*{Acknowledgments}

The authors thank the anonymous referee for a careful, useful, and enthusiastic review. 

We used Claude $4.6$ Sonnet to help connect Equations~\ref{eq:lucaseq2} and \ref{eq:lucaseq4} to identities in Ballot and Williams~\cite{lucas-book}. The LLM never referenced that book. 
This was the only use  we made of AI. All the words are ours.

\section*{Disclosure}

No conflict of interest has been reported by the authors.

\end{document}